\documentclass[12pt]{article}
\usepackage{latexsym} 
\usepackage{amsmath}
\usepackage{amssymb}
\newtheorem{theorem}{Theorem}[section]
\newtheorem{lemma}[theorem]{Lemma}
\newtheorem{corollary}[theorem]{Corollary}
\newtheorem{proposition}[theorem]{Proposition}
\newtheorem{example}[theorem]{Example}
\newtheorem{remark}[theorem]{Remark}

\newtheorem{hypothesis}[theorem]{Hypothesis}

\def\eps{\varepsilon}

\def\bit{\begin{itemize}}
\def\eit{\end{itemize}}
\reversemarginpar   
\def\bc{\begin{center}}
\def\ec{\end{center}}
\def\bthm{\begin{theorem}}
\def\ethm{\end{theorem}}
\def\bcor{\begin{corollary}}
\def\ecor{\end{corollary}}
\def\bprop{\begin{proposition}}
\def\eprop{\end{proposition}}
\def\blem{\begin{lemma}}
\def\elem{\end{lemma}}
\def\bex{\begin{example} {\rm }
\def\eex{\end{example} }}
\def\brem{\begin{remark}}
\def\erem{\end{remark}}
\def\prf{\noindent{\bf Proof. }}
\def\bdes{\begin{description}}
\def\edes{\end{description}}
\def\ita{\item[(a)]}
\def\itb{\item[(b)]}
\def\itc{\item[(c)]}
\def\itd{\item[(d)]}
\def\ite{\item[(e)]}

\def\iti{\item[(i)]}
\def\itii{\item[(ii)]}
\def\itiii{\item[(iii)]}
\def\itiv{\item[(iv)]}

\def\beq{\begin{equation}}
\def\eeq{\end{equation}}
\def\ben{\begin{enumerate}}
\def\een{\end{enumerate}}
\def\beqar{\begin{eqnarray}}
\def\eeqar{\end{eqnarray}}
\def\beqarr{\begin{eqnarray*}}
\def\eeqarr{\end{eqnarray*}}

\def\RR{{\mathbb R}}  

\def\qed{\hspace{.1in}{$\blacksquare$} \\}

\def\Pr{{\mathbf{P}}}

\title{On Gradient like Properties of   Population games, Learning models and Self Reinforced Processes}
 \author{Michel Benaim \\Institut de Math\'ematiques,\\ University of Neuch\^atel}

\begin{document}

\maketitle
\begin{abstract}
We consider ordinary differential equations on the unit simplex of $\RR^n$ that naturally occur in population games, models of learning and self reinforced random processes. Generalizing and relying on an idea introduced in \cite{DF11},  we provide  conditions ensuring that these dynamics are gradient like  and satisfy a suitable "angle condition". This is used to prove that  omega limit sets  and chain transitive sets (under certain smoothness assumptions) consist of equilibria; and that, in the real analytic case, every trajectory converges toward an equilibrium. In the reversible case, the dynamics are shown to be $C^1$ close to a gradient vector field. Properties of equilibria -with a special emphasis on potential games -  and structural stability questions  are also considered.
\end{abstract}
\tableofcontents
\section{Introduction}
Let $S$ be a finite set, say $S = \{1,\ldots, n\}.$ A {\em rate matrix} over $S$ is a $n \times n$ matrix $L$ such that $L_{ij} \geq 0$ for $i \neq j$ and $\sum_{j} L_{ij} = 0.$
We let $\mathbf{R}(S)$ denote the  space of such matrices.   For $x \in \RR^n$ and $L \in \mathbf{R}(S)$ we let $xL$ denote the vector defined by $(xL)_i = \sum_j x_j L_{ji}.$

\noindent \\
Let $$\Delta = \{x \in \mathbb{R}^n : \:  x_i \geq 0, \; \sum_i x_i = 1\}$$ be the {\em unit simplex of probabilities} over $S.$
In this paper we are interested in  ordinary differential equations on $\Delta$ having the form
\beq
\label{ode}
\frac{dx}{dt} = x L(x) : = F(x)
\eeq
where $L : \Delta \mapsto \mathbf{R}(S)$ is a sufficiently smooth function.
Such dynamics  occur - through a natural averaging procedure- in models of games describing strategic interactions in a large population of players, as well as in certain models of learning and reinforcement. These models are usually derived from qualitative assumptions describing the "microscopic" behavior of anonymous agents, and it is usually believed or assumed that similar qualitative microscopic behaviors should lead to similar global dynamics. However there is no satisfactory general theory supporting this belief.

To be more precise, under the assumption that $L(x)$ is irreducible, there exists a unique "invariant probability" for $L(x),$ $\pi(x) \in \Delta$ characterized by
\beq
\label{eq:piL}
\pi(x) L(x) = 0.
\eeq
Several models corresponding to different rate functions  $x \mapsto L(x)$ have the same invariant probability  function $x \mapsto \pi(x).$ For instance,  to each population game (see section \ref{sec:popgames}) which average ODE is given by (\ref{ode}), there is a canonical way to define a learning process (see section \ref{sec:reinforc}) which average ODE is given by
\beq
\label{ode2} \frac{dx}{dt} = - x + \pi(x) := F_{\pi}(x),
\eeq
but there is no evidence that  the dynamics of (\ref{ode}) and (\ref{ode2}) are related in general.

The purpose of this paper is to provide sufficient conditions on $\pi(x)$ ensuring that (\ref{ode}) has a gradient-like structure.  This heavily relies on an idea introduced in \cite{DF11}   where  it was shown that the relative entropy between $x$ and $\pi(x)$ is a strict Lyapounov function for systems of Gibbs type. We extend this idea to other class of systems beyond systems of Gibbs type, including population games and reinforcement process with imitative dynamics, and investigate further dynamics properties.

To give the flavor of the results presented in this paper, let $\pi : \Delta \mapsto \dot{\Delta}$ be a smooth function mapping $\Delta$ into its relative interior. Let  $\chi_{\pi}$ be the set of vector fields having the form given by (the right hand side of) (\ref{ode}), where for each $x,$ $L(x)$ is irreducible and verifies (\ref{eq:piL}). Note that $\chi_{\pi}$ is a convex set of vector fields on $\Delta$ and that $F_{\pi} \in \chi_{\pi}.$
\paragraph{Theorem A}
For all $F \in \chi_{\pi}.$
\bdes
\iti Equilibria (respectively non degenerate equilibria) of $F$ coincide with equilibria (respectively non degenerate equilibria) of $F_{\pi}.$
\itii In general, global dynamics of $F$ and $F_{\pi}$ are "unrelated." We construct an example for which $F_{\pi}$ is globally asymptotically stable (every trajectory converge toward a linearly stable equilibrium) while every non equilibrium trajectory for $F$ converge to a limit cycle.
\itiii Assume that there exists a $C^k, k \geq 1$ strictly increasing function $s:\RR \mapsto \RR$ such that $x \in \dot{\Delta} \mapsto s
((\frac{x_i}{\pi_i(x)}))_{i \in S}$ is the gradient (or  quasi gradient) of some function $V:\dot{\Delta} \mapsto \RR.$ Then
\bdes
\ita $V$ is a strict Lyapounov function for $F$ ($F$ is gradient-like) and verifies an {\em angle condition},
\itb Omega limit sets and chain-transitive sets of $F$ are equilibria,
\itc In the real analytic case, every solution to (\ref{ode}) converge toward an equilibrium,
\itd In the reversible case, hyperbolic equilibria of $F$ coincide with non degenerate critical points of $V$ and, provided there are finitely many equilibria, $F$ is $C^1$ close to a gradient vector field for a certain Riemannian metric,
\ite The set $\chi_{\pi}$ is not (in general) structurally stable.
\edes
\edes

Section \ref{sec:example} describes a few examples that motivate this work. Section \ref{sec:prelim} contains some preliminary results and the main assumptions.
 Section \ref{sec:unrelated} is devoted to Theorem $A, (ii);$
 Section \ref{sec:main} to  $(iii), (a), (b), (c);$ Sections \ref{sec:equilibria} and \ref{sec:reversibility} to $(iii) (d)$ and Section \ref{sec:struct} to $(iii) (e).$
 Other results and examples are also discussed in these sections. For instance, in Section \ref{sec:equilibriapot}, the local dynamics (dynamics in the neighborhood of equilibria) of mean field systems associated
  to potential games is precisely described in term of Nash equilibria.
\section{Motivating  examples}
\label{sec:example}
Throughout this section we see $S$ as set of {\em pure strategies}.
A Markov matrix  over $S$ is a $n \times n$ matrix $K$ such that $K_{ij} \geq 0$ and $\sum_j K_{ij} = 1.$ We let $\mathbf{M}(S)$ denote the sets of such matrices and  we assume given a lipschitz map $$K : \Delta \mapsto \mathbf{M}(S).$$ For further reference we may call such a map a {\em revision protocol.} This terminology is borrowed from \cite{Sand11}.
\subsection{Population Games}
\label{sec:popgames}
Good references on the subject include \cite{Sand10} and the survey paper \cite{Sand11} from which some of the examples here are borrowed.

Consider a population of $N$ {\em agents,} each of whom chooses a strategy in $S$ at discrete times $k = 1,2,\ldots.$ Depending on the context, an agent can be a player, a set of players, a biological entity, a communication device, etc. The state of the system at time $k \in \mathbb{N}$ is the vector $X^N_k = (X^N_{k,1}, \ldots X^N_{k,n}) \in \Delta$ where $N X^N_{k,i}$ equals the number of agents having strategy $i.$ The system evolves as follows. Assume that at time $k$ the system is in state $X^N_k = x.$ Then an agent is randomly chosen in the population. If the chosen agent is an $i-$strategist,  he/she  switches to strategy $j$ with probability  $K_{ij}(x).$
This makes $(X^N_k)_{k \geq 1}$ a discrete time Markov chain, which transition probabilities are
$$\Pr(X^N_{k+1} = x + \frac{1}{N}(e_j-e_i) | X^N_k = x) = x_i K_{ij}(x)$$
where $(e_1, \ldots, e_n)$ is the standard basis of $\RR^n.$
Let
\beq
\label{odepop} L(x) = - Id + K(x).
\eeq
By standard mean-field approximation (see \cite{Kurtz70}, \cite{BW03} for precise statements, and \cite{Sand10}, \cite{Sand11} for discussions in the context of games), the process $\{(X_k^N) \: : k N \leq T\}$ can be approximated by the solution to
(\ref{ode}) (with $L(x)$ given by (\ref{odepop})) with initial condition $x = X_0^N,$ over the time interval $[0,T].$

\subsubsection*{Revision protocols}
Assume, as it is often the case in economic or biological applications, that the population game is determined by a continuous {\em Payoff-function} $U : \Delta \mapsto \RR^n.$ The quantity $U_i(x)$ represents the payoff (utility, fitness) of an $i-$ strategist when the population state is $x.$
\\

\noindent An {\em attachment-function} is a continuous map $w : \Delta \mapsto \RR^{n^2}_+.$ The weight  $w_{ij}(x)$ can be seen as  an  a priori {\em attachment} of an $i-$strategist for strategy $j.$ It can also encompasses certain constraints on the strategy sets. For instance $w_{ij}(x) =0$ (respectively $w_{ij}(x) << 1$) if a move from $i$ to $j$ is forbidden (respectively costly).
We call the attachment function  {\em imitative} if
\beq
\label{eq:defimit}
w_{ij}(x) = x_j \tilde{w}_{ij}(x)
\eeq
Most, not to say all, revision protocols in population games fall into one of the two next categories:
\bdes
\iti [Sampling]
\beq
\label{eq:sampling}
K_{ij}(x) = \frac{w_{ij}(x) f(U_j(x))}{\sum_k w_{ik}(x) f(U_k(x))}
\eeq
where $f$ is a non negative increasing function.
\itii [Comparison]
\begin{eqnarray}
\label{eq:compar}
  K_{ij}(x) &=& w_{ij}(x) g(U_i(x),U_j(x)) \mbox{ for } i \neq j \\
  K_{ii}(x)  &=& 1 - \sum_{j \neq i} K_{ij}(x) \nonumber
\end{eqnarray}
where  $g(u,v)$ is nonnegative, decreasing in $u$, increasing in $v$ and such that $\sum_{j \neq i} K_{ij}(x) \leq 1.$
\edes



\subsection{Processes with reinforcement and adaptive learning}
\label{sec:reinforc}
Suppose now there is only one single agent in the population. In the context of games, one can imagine that this agent consists of a finite set of players and that $S$ is the cartesian product of the strategy sets of the players. let $X_k \in S$ denote the strategy of this agent at time $k.$
Let $\mu_k \in \Delta$ denote the empirical occupation measure of $(X_k)$ up to time $k.$ That is
$$\mu_{k} = \frac{1}{k} \sum_{j = 1}^k \delta_{X_j}$$ where $\delta : S \mapsto \Delta$ is  defined by $\delta_i = e_i.$
Suppose now that the agent revises her strategies as follows:
$$\Pr(X_{k+1} = j | X_0, \ldots, X_{k-1}, X_k = i) = K_{ij}(\mu_k).$$
The process $(X_k)$ is no longer a Markov process but a  {\em process with reinforcement} (see \cite{P07} for a survey of the literature on the subject). Using tools from stochastic approximation theory, it can be shown (see \cite{B97}) that, under certain irreducibility assumptions, the long term behavior of $(\mu_k)$ can be precisely related (see \cite{B97}, \cite{B99}, \cite{BH96} and the brief discussion preceding corollary \ref{chaintrans}) to the long term behavior of the differential equation on $\Delta$
\beq
\label{odereinf}
\frac{dx}{dt} = -x + \pi(x)
\eeq
 where $\pi(x) \in \Delta$ is  the invariant probability of $K(x).$
 Note that (\ref{odereinf}) can be rewritten as  (\ref{ode}) with $L_{ij}(x) = \delta_{ij} - \pi_j(x).$

\section{Hypotheses, Notation, and Preliminaries}
\label{sec:prelim}
Let $L$ be a rate matrix, as defined in the introduction.  Then $L$ is the infinitesimal generator of a continuous time Markov chains on $S.$ A probability $\pi \in \Delta$ is called {\em invariant} for $L$ if it is invariant for the associated Markov chain, or equivalently $$\pi L = 0.$$ A sufficient condition ensuring that $\pi  \in \Delta$ is invariant is that $L$ is {\em reversible} with respect to $\pi,$ meaning that
$$\pi_i L_{ij} = \pi_{j} L_{ji}.$$
The matrix  is  said {\em irreducible} if for all $i,j \in \{1, \ldots, n\}$
there exist some integer $k$ and a sequence of indices $i = i_1, i_2, \ldots, i_{k-1}, i_k = j$   such that
$L_{i_l,i_{l+1}} > 0$    for $l = 1, \ldots, k-1.$

An irreducible rate matrix  admits a unique invariant probability which can be expressed as a rational function of the coefficients $(L_{ij})$ (see e.g~ Chapter 6 of \cite{FW}).

The relative interior of $\Delta$ is the set $$\dot{\Delta} = \{x \in \Delta : \: \forall i \in S, \, x_i > 0\}.$$
 From now on we  assume given a $C^1$ map\footnote{By this we mean that $L$ is the restriction to $\Delta$ of a $C^1$ map defined in a neighborhood of $\Delta$ in  $aff(\Delta) = \{x \in \RR^n \, : \sum_i x_1 = 1\}$} $L : \Delta \mapsto \mathbf{R}(S)$ satisfying the following assumption:
\begin{hypothesis}[standing assumption]
\label{hyp}
For all $x \in \dot{\Delta}, L(x)$ is irreducible.
\end{hypothesis}
We sometimes assume
\begin{hypothesis}[occasional  assumption]
\label{hyp2}
For all $x \in {\Delta}, \, L(x)$ is irreducible.
\end{hypothesis}
In view of the preceding discussion Hypotheses \ref{hyp} and \ref{hyp2}  imply the following
  \blem
  \label{hyperb}   There exists a $C^1$ map $\pi : \dot{\Delta}  \mapsto \dot{\Delta}$
  such that for all $x \in  \dot{\Delta},$ $y = \pi(x)$ is the unique solution  to the equation $$y L(x) = 0, \; y \in \Delta.$$
  If $L$ is $C^k, C^{\infty}$ or real analytic, the same is true for $\pi.$ Under Hypothesis \ref{hyp2}, $\pi$ is defined on all $\Delta$ and maps $\Delta$ into $\dot{\Delta}.$
  \elem
  We let $F_{\pi}$ denote the map defined as
  \beq
  \label{eq:Fpi}
  F_{\pi}(x) = - x + \pi(x).
  \eeq
  Throughout, it is implicitly assumed that the domain of $F_{\pi}$ is $\dot{\Delta}$ under Hypothesis \ref{hyp} and $\Delta$ under Hypothesis \ref{hyp2}.

We now consider the dynamics induced by (\ref{ode}).
 Without loss of generality, we may  assume that (\ref{ode}) is defined on all $\mathbb{R}^n$ and induces a flow $\Phi = (\Phi_t)$ leaving $\Delta$ positively invariant. Indeed, by convexity of $\Delta,$ the retraction $r: \mathbb{R}^n \mapsto \Delta$
  defined by $r(x) = \mathbf{argmin}_{y \in \Delta} \|x-y\|$ is Lipschitz so that the
 differential equation
 \beq
 \label{ode3}
 \frac{dy}{dt} = y L (r(y))
 \eeq is Lipchitz and sub-linear on all $\mathbb{R}^n.$   By standard results,
 it then induces a flow $\Phi : \mathbb{R} \times \mathbb{R}^n \mapsto \mathbb{R}^n$ where
 $t \mapsto \Phi(t,y) = \Phi_t(y)$ is the unique solution to (\ref{ode3}) with initial condition $y.$
 For all $x \in \Delta$ and $t \geq 0,$ $\Phi_t(x) \in \Delta$ and     the map
 $t \in \mathbb{R}^+ \mapsto \Phi_t(x)$ is  solution to (\ref{ode}).

 In the following we sometime use the notation  $\Phi_t(x) = x(t) = ( x_1(t), \ldots, x_n(t)).$
 The {\em tangent space} of $\Delta$ is the space
 $$T\Delta = \{u \in \mathbb{R}^n \: : \sum_{i = 1}^n u_i = 0\}.$$
  \blem
  \label{basiclemma}
 \bdes
 \iti
 There exists $\alpha \geq 0$ such that for all $x \in \Delta$ $x_i(t) \geq e^{-\alpha t} x_i(0).$  In particular, $\dot{\Delta}$ is positively invariant.
 \itii
 If  for all $x \in \partial \Delta$ $L(x)$ is irreducible, then $\Phi_t(\Delta) \subset \dot{\Delta}$ for all $t > 0$ and the dynamics (\ref{ode}) admits a
 global attractor $$A = \bigcap_{t \geq 0} \Phi_t(\Delta) \subset  \dot{\Delta}.$$
 \edes
 \elem
 \prf
 $(i)$ Let $\alpha = \sup_{x \in \Delta} - L_{ii}(x).$ For all $j \neq i$  and $x \in \Delta$
  $$\dot{x}_i \geq - \alpha x_i + x_j L_{ji}(x).$$
  Hence
  $$x_i(t) \geq e^{-\alpha t} [x_i(0) + \int_0^t e^{\alpha s} x_j(s) L_{ji}(x(s))ds] \geq e^{-\alpha t} x_i(0).$$
  The second inequality is the first statement.
  From the first inequality and the continuity of $L(x(t))$ it follows that $x_i(t) > 0$ for all $t > 0 $    whenever that $x_j L_{ji}(x) > 0.$
Let now $x \in \partial \Delta.$ Assume without loss of generality that $x_1 > 0.$ By irreducibility there exists a sequence
 $1 = i_1, i_2, \ldots i_k = j$ such that $L_{i_l,i_{l+1}}(x) > 0.$ Hence, by continuity, $L_{i_l,i_{l+1}}(x(t)) > 0$ for all $t$ small enough.
 It then follows that $x_j(t) > 0$ for all $t > 0$
\qed
\brem {\rm Assumption \ref{hyp} is not needed in  Lemma \ref{basiclemma}.} \erem
Throughout we let $$\mathbf{Eq}(F) = \{x \in \Delta \: : F(x) = 0\}$$
denote the {\em equilibria} set of $F.$
 Note that in view of the preceding Lemmas
 $$\mathbf{Eq}(F) \cap \dot{\Delta}  = \{x  \in  \dot{\Delta} : \: F_{\pi}(x) = 0\}$$ and, in case $L(x)$ is irreducible for all $x \in \Delta,$ $\mathbf{Eq}(F) \subset \dot{\Delta}.$

An equilibrium $p$ is called {\em non degenerate for $F$} provided the Jacobian matrix $DF(p) : T\Delta \mapsto T\Delta$ is invertible.
\blem
\label{lem:equilibria}
 Let $p \in \mathbf{Eq}(F) \cap \dot{\Delta}.$ Then $p$ is non degenerate for $F$ if and only if it is non degenerate for $F_\pi.$
 \elem
 \prf
Let $L^T(x) : T\Delta \mapsto T\Delta$ be defined by $L^T(x) h = h L(x).$ Then for all $x \in \dot{\Delta} \,  F(x) = x L(x) = (x -\pi(x)) L(x) = L^T(x) (x-\pi(x)).$ Hence at every equilibrium $p \in \dot{\Delta}$ $DF(p) = - L^T(p) (DF_{\pi}(p)).$
By irreducibility,  $L^T(p)$ is invertible (see Lemma \ref{lemonLT} in the appendix).  Thus $DF(p)$ is invertible if and only if $DF_{\pi}(p)$ is invertible.
 \qed
\section{Dynamics of $F$ and $F_{\pi}$ are generally unrelated}
\label{sec:unrelated}
While $F$ and $F_{\pi}$ have the same equilibria, they may have quite different dynamics as shown by the following example.

Suppose $n = 3$  so that $\Delta$ is the unit simplex in $\RR^3.$
Let $G$ be a smooth vector field on $\Delta$  such that:
\bdes
\iti $G$ points inward $\dot{\Delta}$ on $\partial \Delta,$
\itii Every forward trajectory of $G$ converge to $p = (1/3, 1/3, 1/3),$
\itiii  $\jmath DG(p) \jmath^{-1} = \left(
                 \begin{array}{cc}
                   - \eta & -1 \\
                   1 & -\eta \\
                 \end{array}
               \right), \, \eta > 0,$

  where $\jmath  : T \Delta \mapsto \RR^2$ is  defined by  $\jmath(u_1,u_2,u_3) = (u_1,u_2).$
\edes
It is easy to construct such a vector field.

Choose $\eps > 0$ small enough so that $\eps G(x) + x$ lies in $\dot{\Delta}$ for all $x  \in \Delta$ and set $\pi(x) =  \eps G(x) + x.$ Then, $F_{\pi}$ and $G$ have the same orbits.

Let $W$ be a $3 \times 3$ symmetric irreducible matrix with positive off-diagonal entries.
Set $L_{ij}(x) = W_{ij} \pi_j(x)$ for $i \neq j$ and $L_{ii}(x) = - \sum_{j \neq i} L_{ij}(x).$
The matrix $L(x)$ is an irreducible rate matrix, reversible with respect to $\pi(x).$
 It follows from Lemmas \ref{basiclemma} and \ref{lem:equilibria} that $F(x) = x L(x)$ has a global attractor contained in $\dot{\Delta}$ and a unique equilibrium given by $p.$

Furthermore, $$DF(p) = - L(p)^T DF_{\pi}(p) = - \eps L(p)^T DG(p) = -\frac{\eps}{3} W DG(p)$$ where the last equality follows from the definition of $L$ and the fact that $\pi(p) = p.$ To shorten notation, set $b = \frac{\eps}{3} W_{12}, c = \frac{\eps}{3} W_{13}$ and
$d = \frac{\eps}{3} W_{23}.$ Then
$$\jmath DF(p) \jmath^{-1}  = \left(
                                                         \begin{array}{cc}
                                                            (b+ 2c) & c-b \\
                                                           d-b &  (b+ 2d) \\
                                                         \end{array}
                                                       \right) \left(
                                                                 \begin{array}{cc}
                                                                   -\eta & -1 \\
                                                                   1 & -\eta \\
                                                                 \end{array}
                                                               \right)$$
The determinant of this matrix is positive and its trace equals $$(c-d) - 2 \eta (b+c+ d).$$
If one now choose $c > d$ and $\eta$ small enough, the trace  is positive. This makes $p$  linearly unstable. By Poincar\'e-Bendixson theorem, it follows that every forward trajectory distinct from $p$ converges toward a periodic orbit.
\brem {\rm It was pointed out to me by Sylvain Sorin and Josef Hofbauer that this example is reminiscent of the following phenomenon. Consider a population game which  revision protocol takes the form
$$K_{ij}(x) =  \frac{x_j}{R} \max(0,U_j(x)-U_i(x)) \mbox{ for } i \neq j$$ (here $R$ is chosen so that $\sum_{j \neq i} K_{ij}(x) \leq  1$). This is a particular case of imitative pairwise comparison protocol (see equation  (\ref{eq:compar})).

Then, the  mean field ode is the classical replicator dynamics (see example 3.2 in \cite{Sand11}):
\beq
\label{replicator}
\dot{x_i} = x_i (U_i(x) - \sum_{j \in S} x_j U_j(x))
\eeq
Here the rate matrix $L(x)$ is not irreducible and its set of invariant probabilities is easily seen to be the {\em Best Reply} set
$$BR(x) = conv (\{ e_i: \: U_i(x) = \max_{j \in S}  U_j(x) \}).$$
where $conv (A)$ stands for the convex hull of $A.$
The vector field (\ref{eq:Fpi}) is not defined but can be replaced by the differential inclusion
\beq
\label{BR}
\dot{x} \in - x + BR(x).
\eeq
If one assume that  $U_i(x) = \sum_{j} U_{ij} x_j$ with  $U$ the payoff matrix given by a Rock-Paper-Scissors game,
$$U = \left(
       \begin{array}{ccc}
         0 & -1 & 1 \\
         1 & 0 & -1 \\
         -1 & 1 & 0 \\
       \end{array}
     \right);$$
     Then $p = (1/3, 1/3, 1/3)$ is the unique  equilibrium of (\ref{replicator}) in $\dot{\Delta}$ (corresponding to the unique Nash equilibrium of the game) and every solution to (\ref{replicator}) with initial condition in $\dot{\Delta}  \setminus \{p\}$ is a periodic orbit.
    On the other hands, solutions to (\ref{BR}) converge to $p.$
    Phase portraits of these dynamics can be found in (\cite{Sand11}, section 5) and a detailed comparison of the replicator and the best reply dynamics is provided in \cite{HSV}.

}
\erem

\section{Gradient Like Structure}
\label{sec:main}
For $u,v \in \RR^n$  we let  $\langle u, v \rangle = \sum_i u_i v_i.$

A map $h :   \dot{\Delta} \mapsto \mathbb{R}^n,$ is called a {\em gradient} if there exists a $C^1$ map $V : \dot{\Delta} \mapsto \mathbb{R}$ such that for all $x \in \dot{\Delta}$ and $u \in T\Delta$
$$\langle h(x), u \rangle = DV(x).u := \langle \nabla V(x), u \rangle.$$

It is called a {\em quasigradient} or a  $\alpha$-{\em quasigradient} if  $x \mapsto \alpha(x) h(x)$ is a gradient for some continuous map $\alpha : \dot{\Delta} \mapsto \mathbb{R}^*_+.$ That is
\beq
\label{defV}
\alpha(x) \langle h(x), u \rangle = \langle \nabla V(x), u \rangle
\eeq
for all $x \in \dot{\Delta}$ and $u \in T\Delta.$
\brem {\rm If $V$ is the restriction to $\dot{\Delta}$ of a $C^1$ map  $W :\RR^n \mapsto \RR,$ then $\nabla V(x)$ is the orthogonal projection of $\nabla W(x)$ onto $T\Delta.$ That is
$$\nabla V_i(x) = \frac{\partial W}{\partial x_i}(x) - \frac{1}{n} \sum_{j = 1}^n \frac{\partial W}{\partial x_j}(x), i = 1, \ldots n.$$}
\erem
\brem  {\rm A practical condition ensuring that $h$ is a gradient is that
\bdes
\ita $h$ is the restriction to $\dot{\Delta}$ of a $C^1$ map $h : \RR^n \mapsto \RR^n,$
\itb For all $x \in \dot{\Delta}$ and $i,j,k \in \{1, \ldots, n\}$
$$\frac{\partial h_i}{\partial x_j}(x) + \frac{\partial h_j}{\partial x_k}(x) + \frac{\partial h_k}{\partial x_i}(x) =
\frac{\partial h_i}{\partial x_k}(x) + \frac{\partial h_k}{\partial x_j}(x) + \frac{\partial h_j}{\partial x_i}(x).$$
This follows from  (\cite{HS98}, Theorem 19.5.5.)
\edes}
\erem

\paragraph{Notation}  We use the following convenient notation. If $x,y $ are vectors in $\RR^n$ and $s : \RR \mapsto \RR$ we let $x.y \in \RR^n$ (respectively $\frac{x}{y}$ and $s(x)$) be the vector defined by $(xy)_i = x_i y_i$ (respectively $(\frac{x}{y})_i = \frac{x_i}{y_i}, s_i(x) = s(x_i)).$

\subsection{Gradient like structure}
A $C^1$ map $V :  \dot{\Delta} \mapsto \RR$ is called a {\em strict Lyapounov function} for $F$ (or $\Phi)$ if for all $x \in \dot{\Delta}$
$$F(x) \neq 0 \Rightarrow \langle F(x), \nabla V(x) \rangle < 0.$$

\bthm
\label{main}
Let  $s : ]0,\infty[ \mapsto \RR$ be a $C^1$ function with positive derivative and let  $h^s :  \dot{\Delta} \mapsto \RR^n$ be the map defined by $$h^s(x) =  s(\frac{x}{\pi(x)}).$$ Assume that $h^s$ is a $\alpha$-quasigradient. Then
\bdes
\iti The map
$V$ (given by (\ref{defV})) is a strict Lyapounov function for $F$ on $\dot{\Delta};$
\itii The critical points of $V$ coincide with  $\mathbf{Eq}(F) \cap \dot{\Delta};$
\itiii $V$ satisfies the following {\em angle condition}: For every compact set $K \subset \dot{\Delta}$ there exists $c > 0$
such that $$\mid \langle \nabla V(x), F(x) \rangle \mid  \geq c \parallel \nabla V(x) \parallel \parallel F(x) \parallel $$ for all $x \in K.$
\edes
\ethm
\brem{\bf [Gibbs systems]}

 \label{rem:DF}
 {\rm If $\pi(x)$ is a {\em Gibbs measure},
 \beq
 \label{Gibbs}
 \pi_{\beta,i}(x) = \frac{ e^{- (U_i^0 +  \beta \sum_j U_{i j} x_j)}}{Z(x)}
  \eeq
  where $U = (U_{ij})$ is a  symmetric matrix, $\beta \geq  0,$ and $$Z(x) = \sum_j e^{ - (U_j^0 +  \beta \sum_k U_{i k} x_k)},$$ parts
$(i)$ and $(ii)$  of Theorem \ref{main} have been  proved in \cite{DF11}, Theorems 5.3 and 5.5. Here $s(t) = \log(t)$ and
\beq
\label{freeenergy}
V(x) = \sum_i x_i \log(x_i) + \sum_j U_j^0 x_j +  \frac{\beta}{2} \sum_{ij} U_{ij}x_i x_j.
\eeq}
\erem
\subsubsection*{Proof of theorem \ref{main}}

Part $(i)$  relies on the following Lemma.
\blem
\label{poincare1}
Let $L$ be an irreducible transition matrix with invariant probability $\pi.$ Let $x \in \Delta, f_i = \frac{x_i}{\pi_i}, \, s(f)_i = s(f_i)$ and $c_f = \inf_{i} s'(f_i) > 0.$ Then there exists $\lambda(L) > 0$  depending continuously  on $L$ such that $$\langle x L, s(f) \rangle \leq -c_f \lambda(L) Var_{\pi}(f)$$
where $Var_{\pi}(f) = \sum_i (f_i - 1)^2 \pi_i = \sum_i \frac{(x_i -\pi_i)^2}{\pi_i}.$
\elem
The proof of this lemma uses  elementary convexity arguments and classical tools from Markov chain theory. It is proved in appendix.
Applying this lemma with $L = L(x)$ and  $\pi = \pi(x)$ gives $$\langle F(x), \nabla V(x) \rangle < 0$$ unless $x = \pi(x).$

$(ii)$ The set $\mathbf{Eq}(F) \cap \dot{\Delta}$ coincides with fixed points of $\pi$ in $\dot{\Delta}.$
Let  $x \in \dot{\Delta}.$ $\nabla V(x) = 0 \Leftrightarrow h^s(x) \in \RR 1$ where $1$ is the vector which components are all
 equal to $1.$ The function $s$ being injective this is equivalent to $\frac{x_i}{\pi_i(x)} = \frac{x_j}{\pi_j(x)}$ for all $i,j.$ That is
 $x = \pi(x).$

 $(iii)$ Let $K \subset \dot{\Delta}.$ By   Lemma \ref{poincare1} (applied with $L = L(x)$ and $\pi = \pi(x)$) and   continuity of the maps  involved, there exists $c > 0$ depending on $K$ such that
 $$| \langle \nabla V(x), F(x) \rangle | \geq c \sum_i (x_i - \pi_i(x))^2 = c \|x-\pi(x)\|^2.$$
 To prove the angle condition it then suffices to show that both $\|F(x)\|$ and $\|\nabla V(x)\|$ are bounded by some constant times $\|x-\pi(x)\|.$
 Now, $F(x) = x L(x) = x L (x) - \pi(x) L(x)$ so that $$\|F(x)\| \leq c_1 \|x-\pi(x)\|$$ with $c_1 = \sup_{x \in \Delta} \|L(x)\|.$

By Lipschitz continuity of $s$ and compactness, there exist $c_2, c_3 > 0$ depending on $K$ such that
 $$|s(\frac{x_i}{\pi_i(x)}) - s(1)| \leq c_2 |\frac{x_i}{\pi_i(x)} - 1| \leq c_3 |x_i - \pi_i(x)|.$$
 Thus, for all  $u \in T\Delta$ such that $\|u\| = 1$
 $$\langle h^s(x) , u \rangle = \langle h^s(x) - s(1)1, u \rangle \leq \|h^s(x) - s(1).1\| \leq c_3 \|x-\pi(x)\|.$$
 This implies that $\|\nabla V(x)\| \leq c_3 \|x-\pi(x)\|$ and concludes the proof. \qed

The following result proves to be useful for certain dynamics leaving invariant the boundary of the simplex. Such dynamics occur in population games using imitative protocols (see equation \ref{eq:defimit}) as well as in certain models of vertex reinforcement (see example \ref{ex:VRRW} below).

For $x \in \Delta$ let $Supp(x) = \{x \in \Delta \: : x_i > 0\}.$
\bprop
\label{main2}
 Assume that assumptions of Theorem \ref{main} hold. Assume furthermore that
\bdes
\ita For all $x \in \Delta$
$$x_i = 0 \Rightarrow L_{ji}(x) = 0$$ and the reduced rate matrix $[L_{ij}(x)]_{i,j \in Supp(x)}$ is irreducible
\itb
 The maps $V : \dot{\Delta} \mapsto \RR^n$  and $\alpha : \dot{\Delta} \mapsto \RR^*_+$ (given by equation \ref{defV}) extend to  $C^1$ (respectively continuous) maps $V : \Delta \mapsto \RR^n$ and $\alpha : \Delta \mapsto \RR^*_+.$

\edes
Then $V$ is strict Lyapounov function for $F$ on $\Delta.$
\eprop
\prf Let $T \Delta(x) = \{u \in T \Delta \: : u_i = 0 \mbox{ for } i \not \in Supp(x)\}.$ By assumption $(a)$ the map $x \mapsto \pi(x)$ is defined for all $x \in \Delta$ continuous in $x$ and $\pi_i(x) = 0 \Leftrightarrow x_i = 0.$

Therefore, using assumption $(b)$, the equation
$$ \forall x \in \dot{\Delta}, \forall u \in T\Delta \: \alpha(x) \langle h^s(x), u \rangle = \langle \nabla V(x), u \rangle$$
extends to
$$\forall x \in \Delta, \forall u \in T\Delta(x) \: \sum_{i \in Supp(x)} s(\frac{x_i}{\pi_i(x)}) u_i = \langle \nabla V(x), u \rangle$$
Thus $$\sum_{i \in Supp(x)} s(\frac{x_i}{\pi_i(x)}) (xL(x))_i = \langle \nabla V(x), F(x) \rangle$$ for all $x \in \Delta.$
By Lemma \ref{poincare1} the left hand side is nonpositive and zero if and only if $x_i = \pi_i(x)$ for all $i \in Supp(x).$  \qed
\brem {\rm Note that under the assumptions of Proposition \ref{main2}, the angle inequality of Theorem \ref{main} doesn't hold on the boundary of the simplex }\erem
\bex
\label{ex:gener}
{\rm
Let $W : \RR^n \mapsto \RR$ be a $C^k$ map, $k \geq 1$.
Suppose that for all $x \in \dot{\Delta}$
$$\pi_i(x) =  \frac{f_i(x_i) \psi(\frac{\partial W}{\partial x_i}(x))}{\sum_{j = 1}^n f_j(x_j) \psi(\frac{\partial W}{\partial x_j}(x))} $$
Then, Theorem \ref{main} applies in the following cases:
\bdes
\item{{\bf Case 1}}
$\psi(u) = e^{- \beta u}$ with $\beta \geq 0,$ and  $f_i(t) > 0$ for all $t > 0.$   It suffices to choose $s(t) = \log(t)$ and
\beq
\label{eq:lyapgibbs} V(x) = \sum_{i = 1}^n x_i \log(x_i) - \sum_{i = 1}^n \int_1^{x_i} \log(f_i(u)) du + \beta W(x)
\eeq
Then $h^s$ is the gradient of $V.$
\item{{\bf Case 2}} $\psi(u) = u^{\beta}, \beta > 0, f_i(t) = t$ and
$\frac{\partial W}{\partial x_i} > 0$ on $\{x \in \Delta : \: x_i > 0 \}.$ It suffices to choose  $s(t) = - t^{- 1/\beta}$ and
$$V(x) = - W(x).$$
Then $h^s$ is the $\alpha$-quasigradient of $V$ with $$\alpha(x) = [\sum_j x_j (\frac{\partial W}{\partial x_j})^{\beta})]^{-1/\beta}.$$
\edes
}
\eex
\bex{\bf [Potential Games]}
{\rm
\label{sec:potential}

Examples of applications of   Example \ref{ex:gener}, case 1, are given by {\em Potential Games} (see \cite{Sand10} for an comprehensive presentation and motivating examples).
We use the notation of section \ref{sec:example}. A {\em Potential Game} is a game for which the payoff function is such that for all $x \in \Delta$
 $$U_i(x) = - \frac{\partial W}{\partial x_i}(x), i = 1 \ldots n$$ 
 Consider a population game with a revision protocol given by (\ref{eq:compar}).
 Suppose that  the attachment matrix takes the form
 $$w_{ij}(x) = f_j(x_j) \tilde{w}_{ij}(x)$$
 with $\tilde{w}$ irreducible and symmetric. Let $\beta \geq 0.$ Assume furthermore that $g(u,v)$ takes one of the following form:
 \\
 {\em Pairwise comparison}
 $$g(u,v) = \frac{e^{\beta (v-u)}}{1+e^{\beta (v-u)}} \mbox{ or } g(u,v) =  \min(1,e^{\beta (v-u)}),$$
 {\em Imitation driven by dissatisfaction}
 $$g(u,v) =
 e^{- \beta u},$$
  {\em Imitation of success}
$$g(u,v) =  e^{ \beta v}.$$
 In all these situations, $K(x)$, hence $L(x)$ is reversible with respect to $\pi_{\beta}(x)$ with
 $$\pi_{\beta,i}(x) = \frac{f_i(x_i) e^{-\beta \frac{\partial W}{\partial x_i}(x)}}{\sum_{j} f_j(x_j) e^{-\beta \frac{\partial W}{\partial x_j}(x)}}.$$
 Theorem \ref{main} applies with $V$ given by (\ref{eq:lyapgibbs}).
 \brem{\bf [Gibbs systems, 2]}
 \label{rem:DF2}{\rm A particular case of potential games is obtained with
 $W(x) = \frac{1}{2} \sum_{ij} U_{ij} x_i x_j$ with $U = (U_{ij})$ symmetric, and $f_i(x) = e^{-U^0_i}.$
  Here payoffs are linear in $x:$
  $$U_i(x) = - \sum_j U_{ij} x_j$$ and we retrieve  the situation considered in \cite{DF11}. See Remark \ref{rem:DF}.}
 \erem
 }
 \eex

\bex{\bf [Vertex reinforcement]}
\label{ex:VRRW}
{\rm
Let $K$ be the revision protocol defined by
 $$K_{ij}(x) = \frac{A_{ij} x_j^{\gamma}}{\sum_k A_{ik} x_k^{\gamma}}$$
 where $A$ is a  matrix with positive entries and $\gamma \geq 1.$ For population games (see section \ref{sec:popgames}) this gives a  simple model of imitation: an agent of type $i,$ when chosen,  switches to $j$ with a probability  proportional to the $( \mbox{{\em number of agents of type } }j)^{\gamma}.$
 For processes with reinforcement (as defined in section \ref{sec:reinforc})  the probability to jump from $i$ to $j$ at time $n$ is proportional to
 $(\mbox{{\em
  the time spent in }} $j$ \mbox{ {\em up to time} } $n$)^\gamma.$ This later model  called a {\em vertex reinforced random walks}  was introduced by Diaconis and first analyzed in Pemantle \cite{P92} (see also \cite{B97} and \cite{BRS} for more references on the subject).

  When $A$ is symmetric, $K(x)$ is reversible with respect to
 \beq
 \label{pireinf}
 \pi_i(x) = \frac{x_i^{\gamma} \sum_k A_{ik} x_k^{\gamma}}{\sum_{ij} A_{ij} x_i^{\gamma} x_j^{\gamma}} =
 \frac{ x_i \frac{\partial  W}{\partial x_i}}{\sum_j x_j \frac{\partial  W}{\partial x_j}}
 \eeq
 with
 \beq
 \label{eq:defVreinf}
 W(x)  =  \sum_{i,j} A_{ij} x_i^{\gamma} x_j^{\gamma}
 \eeq
  We are then in the situation covered by  Example \ref{ex:gener}, case 2, with
 $\psi(u) = u, f_i(t) = t, s(t) = -\frac{1}{t}$ and $V = - W.$

 Both Theorem \ref{main} and Proposition  \ref{main2} apply.
}
 \eex
 \bex{\bf [Interacting Urn processes]}
 {\rm Closely related to vertex reinforced random walks are models of {\em interacting urns}  (see
 \cite{BBCL}, \cite{CL13}, \cite{HofHolmes14}). For these models $\pi_i(x) = x_i \frac{\partial W}{\partial x_i}$ for some smooth function $W.$ This is a particular case of Example \ref{ex:gener}, case 2.}
 \eex

\subsection{Limit sets and chain transitive sets}
\label{sec:limset}
Using Lasalle's invariance principle we deduce the following consequences from Theorem \ref{main}.
\bcor
\label{lassale} Assume that assumptions of Theorem \ref{main} hold. Then every  omega limit set
 of $\Phi$  contained in $\dot{\Delta}$
 is a connected subset of $\mathbf{Eq}(F) \cap \dot{\Delta}.$
\ecor
Combining this results with Lemma \ref{basiclemma} (ii) and Proposition \ref{main2} gives
\bcor Assume that one of the following condition hold:
\bdes \ita  Assumptions of Theorem \ref{main} and Hypothesis \ref{hyp2} or;
\itb Assumptions of Proposition \ref{main2}.
\edes
Then every  omega limit set
 of $\Phi$
 is a connected subset of $\mathbf{Eq}(F).$
\ecor
A set $L$ is called {\em attractor free} or {\em internally chain transitive} provided it is compact, invariant and $\Phi|L$ has no proper attractor.
For reinforced random processes like the ones defined in section \ref{sec:reinforc}, limit sets of $(\mu_n)$ are, under suitable assumptions, attractor free sets of the associated mean field equation (\ref{odereinf}) (see \cite{B97}). More generally attractor free sets are limit sets of {\em asymptotic pseudo trajectories} (see \cite{BH96}). It is then useful to characterize such sets.  Note however that the existence of a strict Lyapounov function, doesn't ensure in general, that internally chain transitive sets consist of equilibria (see e.g~  Remark 6.3 in \cite{B99}).
\bcor
\label{chaintrans}
Assume that assumptions of Theorem \ref{main} hold and that $h^s$ is $C^{k}$ for some $k \geq n-2 = \dim( T \Delta) - 1.$ Then every internally chain transitive set of $\Phi$  contained in $\dot{\Delta}$
 is a connected subset of $\mathbf{Eq}(F) \cap \dot{\Delta}.$
 If we furthermore assume that  $L(x)$ is irreducible for all $x \in \Delta,$ then every internally chain transitive set of $\Phi$
 is a connected subset of $\mathbf{Eq}(F)$
\ecor
\prf Let $C = \mathbf{Eq}(F) \cap \dot{ \Delta }$ and $A \subset \dot{\Delta}$ an attractor free set. By Theorem \ref{main}, $C$ coincide with critical points of $V.$ By the assumption $V$ is $C^{k+1}$ so that by  Sard's theorem (see \cite{Hirsch76}), $V(C)$ has empty interior. It  follows (see e.g Proposition 6.4 in \cite{B99}) that $A \subset C.$
\qed
\subsection{Convergence toward one equilibrium}
In case equilibria are isolated, Corollary \ref{lassale} implies that every trajectory bounded away from the boundary converge to an equilibrium and that every trajectory converges in case $L(x)$ is irreducible for all $x.$
However, when equilibria are degenerate, the gradient-like property is not sufficient to ensures convergence. There are known examples of smooth gradient systems which omega limit sets are a continuum of equilibria (see \cite{Palis}).
However, in the real analytic case, gradient like systems which verify an angle condition are known to converge
\bthm
Suppose  that assumptions of Theorem \ref{main} hold and that $V$ is real analytic. Then every omega limit set meeting $\dot{\Delta}$ reduces to a single point.
\ethm
\prf Let $p$ be an omega limit point. If $V$ is real analytic, it satisfies a {\em Lojasiewicz inequality}  at $p$ in the sense that there exist $0 < \eta \leq 1/2, \beta > 0$ and a neighborhood $U(p)$ of $p$ such that
$$|V(x) - V(p) |^{1-\eta} \leq \beta \|\nabla V(x)\|$$ for all $x$ in a $U(p)$. Such an inequality called a "gradient inequality"
was proved by Lojasiewicz \cite{Loja}
 and used (by Lojasiewicz again) to show that  bounded solutions of real analytic gradient vector fields have finite length, hence converge. When the dynamics is not a gradient, but only gradient like with $V$ as a strict Lyapounov function, the same results holds provided that $V$ satisfies an angle condition:
$$\langle \nabla V(x), F(x) \rangle \geq c \|F(x)\|\|\nabla V(x)\|$$ for all $x \in U(p).$  This is proved in \cite{chill09} (see also  \cite{Merlet13}, Theorem 7). \qed
\bex {\bf [Gibbs systems, 3]}
 {\rm If $\pi$ is given by (\ref{Gibbs}) with $U$ symmetric, $V$ given by (\ref{freeenergy}) is real analytic
 so that every solution to (\ref{ode}) converges toward an equilibrium.}
 \eex
\section{Equilibria}
\label{sec:equilibria}
Recall that point $p \in \mathbf{Eq}(F)$ is called {\em non degenerate} if the jacobian matrix $DF(p):T\Delta \mapsto T\Delta$ is invertible. It is called {\em hyperbolic} if eigenvalues of $DF(p)$ have non zero real parts. If $p$ is hyperbolic,  $T\Delta$ admits a splitting $$T\Delta = E^u_p \oplus E^s_p$$ invariant under $DF(p)$ such that the eigenvalues of  $DF(p)|_{E^s_p}$ (respectively   $DF(p)|_{E^u_p}$) have negative (respectively positive) real parts.

Point $p \in \mathbf{Crit}(V) = \{x \in \dot{\Delta} \: : \nabla(V)(x) = 0\}$ is called {\em non-degenerate} if $Hess(V)(p)$  the {\em Hessian} or $V$ at $p$  has full rank. In a suitable coordinate systems
$Hess(V)(p) (u,u)  = \sum_{i = 1}^{n_+} u_i^2 - \sum_{j = 1}^{n_-} u_j^2$ with $n_+ + n_{-} = dim(T\Delta) = n-1.$
 The number $n_-$ is called the {\em index} of $p$ (with respect to $V$) and is written $\mathbf{Ind}(p,V).$
\bprop
\label{prop:equilibria}
Assume that assumptions of Theorem \ref{main} hold. Let $p \in \mathbf{Eq}(F) \cap \dot{\Delta}.$ Then
\bdes
\iti Point $p$ is non degenerate if and only if it is a non degenerate critical point of $V.$
\itii If furthermore $L$  is $C^2$ and $p$ is hyperbolic,  $$\dim (E^u_p) = \mathbf{Ind}(p,V).$$
\edes
\eprop
\prf
From Lemma \ref{lem:equilibria}, $p$ is non degenerate if and only if $DF_{\pi}(p)$ is invertible and (see the proof of Lemma \ref{lem:equilibria})
\beq
\label{dfp}
DF(p) = -L^T(p) DF_{\pi}(p)
 \eeq
Now, a direction computation (details are left to the reader)  of the Hessian of $V$ at $x$
leads to
$$\langle Hess(V)(x) u,v \rangle = \alpha(x) \langle s'(\frac{x}{\pi(x)}) (u- \frac{x}{\pi(x)} D\pi(x) u, v \rangle_{1/\pi(x)}$$
where   $\langle u,v \rangle_{1/{\pi}}$ stands for  $\sum_i u_i v_i \frac{1}{\pi_i}.$ Since $p = \pi(p)$
\beq
\label{defHess}
\langle Hess(V)(p) u,v \rangle = \alpha(p) s'(1) \langle (I - D\pi(p)) u, v \rangle_{1/p}
\eeq for all $u,v \in T\Delta,$
This proves that $Hess V(p)$ is non degenerate if and only if  $(I - D\pi(p)) = - DF_{\pi}(p)$ is non degenerate and concludes the proof of the first part.

We now prove the second part. By the stable manifold theorem, there exists a  (local)  $C^2$ manifold  $W^s_p$ tangent to $E^s_p$ at $p$ positively invariant under $\Phi$ and such that for all $x \in W^s_p$ $\lim_{t  \rightarrow  \infty} \Phi_t(x) = p.$ Clearly $p$ is a global minimum of $V$ restricted to $W^s.$ For otherwise there would exists $x \in W^s_p$ such that  $$V(p) > V(x) > \lim_{t \rightarrow \infty} V(\Phi_t(x)) =  V(p).$$
Since $p$ is also a critical point $\nabla V(p) = 0.$ Let $u \in E^s_p$ and let $\gamma : ]-1,1[ \mapsto W^s_p$ be a $C^2$ path with $\gamma(0) = p, \dot{\gamma}(0)  = u.$ Set $h(t) = V(\gamma(t)).$ Then $\dot{h}(0) = 0$ (because $p$ is a critical point of $V$) and $h''(0) = \langle Hess V(p) u, u \rangle $ is non negative because   $h(t) \geq h(0).$

On the other hand, by the spectral decomposition of $Hess V(p)$ we can write $T \Delta = E_V^s \oplus E_V^u$ with  $ \langle Hess V(p) u , u \rangle > 0$ (respectively $< 0$) for all $u \in E_V^s \setminus \{0\}$ (respectively $E_V^u \setminus \{0\}$).
Thus, $E_p^s \cap E_V^u = \{0\}$ and, consequently,  $\dim(E_p^s) + \dim (E_V^u) \leq \dim(T\Delta).$ Similarly $\dim(E_p^u) + \dim (E_V^s) \leq \dim(T\Delta).$ This proves that $\dim(E^u_p) = \dim (E_V^u) = \mathbf{Ind}(p,V).$
 \qed

\brem
\label{hartman}
{\rm This later proposition shows that in the neighborhood of an hyperbolic equilibrium $p$, $\dot{x} = F(x)$ and $\dot{x} = -\nabla V(x)$ are topologically conjugate. Indeed, part $(ii)$ of the proposition implies that the linear flows $\{e^{t DF(p)}\}$ and $\{e^{t Hess(V)(p)}\}$  are topologically conjugate (see e.g~ Theorem 7.1 in \cite{R1}), and by Hartman-Grobman Theorem (see again \cite{R1}), nonlinear flows are locally conjugate to their linear parts in the neighborhood of hyperbolic equilibria. However, note that while eigenvalues of $Hess(V)(p)$ are reals there is no evidence that the same is true for $DF(p)$ in general. The next proposition proves that this is the case when  $L(x)$ is reversible with respect to $\pi(x).$ }\erem

\bprop
\label{reversequil}
Let
$p \in \mathbf{Eq}(F) \cap \dot{\Delta}.$ Assume that assumptions of Theorem \ref{main} hold  and that
 $L(p)$ is reversible with respect to $\pi(p) = p.
$ Then there exists a positive definite bilinear  form $g_0(p)$ on  $T\Delta$  such that for all $u,v \in T\Delta$
$$g_0(p)(DF(p) u, v) = - \langle Hess(V)(p) u, v \rangle$$
In particular \bdes
\iti $DF(p)$ has real eigenvalues,
\itii $p$ is hyperbolic for $F$ if and only if it is a non degenerate critical point of $V.$
\edes
\eprop
\prf
Let $p \in  \mathbf{Eq}(F) \cap \dot{\Delta}.$
 Set $L = L(p)$ and recall that   $L^T : T \Delta \mapsto T \Delta$ is defined by $ L^T h =  h L.$
Then, by Lemma  \ref{lemonLT} in the appendix,
$-L^T$ is a definite positive operator for  the scalar product on $T \Delta$ defined by
$\langle u, v \rangle_{1/p} = \sum_i u_i v_i \frac{1}{p_i}.$ Define now $g_0(p)$ by
\beq
\label{eq:metric} g_0(p)(u,v) = - \langle (L^T)^{-1} u, v \rangle_{\frac{1}{p}}.
\eeq
Using (\ref{dfp}) and (\ref{defHess}) it comes that for all $u,v \in T\Delta$
$$g_0(p)(DF(p) u, v) = - \langle (I - D\pi(p)) u, v \rangle_{1/p}$$
 $$= - [\alpha(p) s'(1)]^{-1} \langle Hess(V)(p) u,v \rangle.$$ Replacing $g_0(p)$ by $\alpha(p) s'(1) g_0(p)$ proves the result.\qed

  A useful consequence of this later proposition is that it is usually much easier to verify non degeneracy of  equilibria rather than hyperbolicity. Here is an illustration:
 \bex{\bf [Gibbs systems, 4]}
 \label{gibbs}{\rm
 Consider the symmetric Gibbs model analyzed in \cite{DF11} (see remark \ref{rem:DF} and example \ref{rem:DF2}). We suppose that the symmetric matrix $U = (U_{ij})$ is given and we treat $U^0 = (U^0_i)_{i \in S}$ and $\beta$ as parameters.  Let ${\Xi}_{rev}(U^0)$ denote the set of maps
 $$\RR^+ \times \Delta \mapsto T\Delta,$$
 $$(\beta,x) \mapsto F_{\beta}(x) = x L_{\beta}(x)$$
 such that $L_{\beta}$  verifies assumption \ref{hyp2}, is $C^1$ in $x,$ and $L_{\beta}(x)$ is reversible with respect to $\pi_{\beta}(x)$ where
 $$\pi_{\beta,i}(x) = \frac{e^{- U^0_i -    \beta \sum_{j} U_{ij} x_j }}{Z(x)}.$$
 \bprop There exists  an open and dense set ${\cal G}^0 \subset \RR^n$ such that for all $U^0 \in {\cal G}^0$ and $F \in \Xi_{rev}(U^0)$
 \bdes
 \iti The set $\{(x,\beta) \in \Delta \times \RR^+_: : \: F_{\beta}(x) = 0 \}$ is a $C^{\infty}$ one dimensional submanifold,
 \itii There exists an open dense set ${\cal B}^0 \subset  \RR^+ $ containing $0$ such that for all $\beta \in {\cal B}^0$ equilibria of $F_{\beta}$ are hyperbolic.
 \edes
 \eprop
 \prf
 Let $H : \dot{\Delta} \times \RR^n \times \RR^+ \mapsto T\Delta$ be defined by
 $H(x,U^0,\beta) = \nabla V_{U^0,\beta}(x)$ where $V_{U^0,\beta}$ is given by (\ref{freeenergy}).
  Since $\frac{\partial H}{\partial U^0}(x,U^0,\beta)$  is the identity map, $H$ is a submersion. Hence, by Thom's parametrized transversality Theorem (see \cite{Hirsch76}, Chapter 3), there exists an open and dense set ${\cal G}^0 \in \RR^n$ such that for all $U^0 \in {\cal G}^0, (x,\beta) \mapsto H(x,U^0,\beta)$ is a submersion. This proves $(i).$ By the same theorem, for  all $\beta \in {\cal B}^0$ with ${\cal B}^0$ open and dense in $\RR^+, x \mapsto H(x,U^0,\beta)$ is a submersion, meaning that critical points of $V_{U^0,\beta}$ are nondegenerate. By Proposition \ref{reversequil}, equilibria of $F_{\beta}$ are hyperbolic.
 \qed

 \brem {\rm Other genericity results can be proved, if one fix $U^0$ or $\beta$ and treat $U$ as a parameter.
  Compare to the proof of Theorem 2.10 in \cite{BR05} in an infinite dimensional setting. }
  \erem
  }
 \eex
\subsection{Equilibria of Potential Games}
\label{sec:equilibriapot}
Consider a population game with $C^1$ payoff function $U: \Delta \mapsto \RR^n.$ Recall that the game is called a potential game, provided $U_i(x) = - \frac{\partial W}{\partial x_i}(x)$ for all $x \in \Delta$ and some {\em potential} $W : \RR^n \mapsto \RR.$

Point $x^* \in \Delta$ is called a {\em Nash equilibrium} of $U$ if, given the population state $x^*$, every agent has interest to play the mixed strategy $x^*.$ That is
\beq
\label{eq:nash}
\forall i \in \{1, \ldots, n \} \:   U_i(x^*)  \leq \langle U(x^*), x^* \rangle
\eeq
Let
$$Supp(x^*) = \{ i \in \{1, \ldots, n\}: \: x^*_i > 0 \}.$$ It follows from (\ref{eq:nash}) that
$$ \forall i \in Supp(x^*) \:   U_i(x^*)  = \langle U(x^*), x^* \rangle.$$

We let ${\mathbf {NE}}(U)$  denote the set of Nash equilibria of $U.$
For all $\beta \geq 0$ and $x  \in \Delta$ we let $\pi_{\beta}(x) \in \Delta$ be defined as
\beq
\label{defpibeta}
\pi_{\beta,i}(x) = \frac{e^{\beta U_i(x)}}{\sum_j e^{\beta U_j(x)}}, \, i = 1, \ldots, n
\eeq
and
we let
$\chi(\beta, U)$ (respectively, $\chi_{rev}(\beta, U))$ denote the set of all vector fields having the form given by (\ref{ode}) where $L(x)$ is $C^1$ in $x$, irreducible and admits $\pi_{\beta}(x)$ as invariant (respectively reversible) probability. Recall (see equation (\ref{eq:Fpi})) that  $$F_{\pi_{\beta}} = - Id + \pi_{\beta}.$$
Our aim here is to describe $\mathbf{Eq}(F)$ for  $F \in \chi(\beta, U)$ in term of $\mathbf{NE}(U)$ for large $\beta,$ with a particular emphasis on potential games.  Some of the results here are similar to the results obtained in \cite{BH06} for $n \times 2$ pseudo games.

\bprop
\label{prop:equiNash}
Let ${\cal N}$ be a neighborhood of ${\mathbf{NE}}(U).$  There exists $\beta_0 \geq 0$ such that for all $\beta \geq \beta_0$ and $F \in \chi(\beta,U)$ $${\mathbf{Eq}}(F)  \subset  {\cal N}$$
\eprop
\prf Equilibria of $F$  coincide with equilibria of  $F_{\pi_{\beta}}.$ Let $x(\beta)$ be such an equilibrium. Then  for all $i,j$
$$\frac{\log(x_i(\beta)) - \log(x_j(\beta))}{\beta} = U_i(x(\beta)) - U_j(x(\beta)).$$
Thus for every limit point $x^* = \lim_{\beta_k \rightarrow \infty} x(\beta_k)$ it follows that
$$U_i(x^*) = U_j(x^*)$$ if $i,j \in Supp(x^*)$ and
$$U_i(x^*) \geq U_j(x^*)$$ if $i \not \in Supp(x^*)$ and $j \in Supp(x^*).$
\qed
\brem {\rm Note that Proposition \ref{prop:equiNash} only requires the continuity of $U.$ } \erem

We shall now prove some converse results.

A Nash equilibrium $x^*$ is called {\em pure} if $Supp(x^*)$ has cardinal $1$ and {\em mixed} otherwise. It is called {\em strict} if inequality (\ref{eq:nash}) is strict for all $i \not \in Supp(x^*).$

\bthm
\label{thm:equipure}
Let $x^*$ be a pure strict Nash equilibrium and  ${\cal N}$ a (sufficiently small)  neighborhood of $x^*.$  Then, there exists $\beta_0 > 0$ such that for all $\beta \geq \beta_0$ and   $F \in \chi(\beta,U)$
\bdes
\iti $\mathbf{Eq}(F) \cap {\cal N} = \{x^*_{\beta}\}$
\itii Equilibrium $x^*_{\beta}$ is linearly stable for $F_{\pi_{\beta}}.$
\itiii
Assume furthermore that the game is a potential game. Then $x^*_{\beta}$ is linearly stable for $F$ under one of the following conditions: \bdes \ita $L$ (hence $F$) is $C^2$ and $x^*_{\beta}$ is hyperbolic for $F,$  or
\itb $F \in \chi_{rev}(\beta, U).$ \edes
\edes
\ethm
\prf
Suppose without loss of generality that $x^*_1 = 1$ and $x^*_i = 0$ for $i \neq 1.$

 Set $R_{ij} = U_{j} - U_{i}.$
By assumption and continuity, there exists $\delta > 0, \alpha > 0$ such that for all $x \in B(x^*,\alpha) = \{ x \in \Delta \: : \|x-x^*\| \leq \alpha\},$
$$R_{i1}(x) \geq \delta  \mbox{ for } i > 1;$$
 $$\|R_{ij}(x)\| \geq \delta \mbox{ if } R_{ij}(x^*) \neq 0$$ and $$\|R_{ij}(x)\| \leq \delta \mbox{ if } R_{ij}(x^*) = 0.$$

 Thus $$1 \geq \pi_{\beta,1}(x) = (1 + \sum_{i > 1} e^{- \beta R_{i1}(x)})^{-1} \geq (1 + (n-1) e^{-\beta \delta})^{-1}.$$
 This implies that $\pi_{\beta}$ maps $B(x^*, \alpha)$ into itself for $\beta$ large enough. By Brouwer's Theorem, it then admits a fixed point $x^*_{\beta}.$ To prove uniqueness and assertion $(ii)$ it suffices to prove that $\pi_{\beta}$ restricted to $B(x^*,\alpha)$ is a contraction.
From the expression $\pi_{\beta,i} = (\sum_{j} e^{\beta R_{ij}})^{-1},$ we get
 $$\frac{\partial \pi_{\beta,i}}{\partial x_m} = - \sum_{j} [\beta e^{\beta R_{ij}}(\sum_{k} e^{\beta R_{ik}})^{-2} \frac{\partial R_{ij}}{\partial x_m}] := \sum_{j} D_{ij} = \sum_{j \neq i} D_{ij}.$$
Let  $j \neq i.$ If $R_{ij}(x^*) \neq 0$
 $$|D_{ij}| \leq \beta e^{\beta R_{ij}} (1 + e^{\beta R_{ij}})^{-2} \leq \beta \min(e^{\beta R_{ij}}, e^{-\beta R_{ij}}) \leq \beta e^{-\beta \delta}.$$
 If $R_{ij}(x^*) = 0.$ Then $i \neq 1$ and
 $$|D_{ij}| \leq \beta e^{\beta R_{ij}} (e^{\beta R_{i1}})^{-2} = \beta e^{\beta (R_{ij} - 2 R_{i1})} \leq \beta e^{-\beta \delta}$$
These inequalities show that $\|D\pi_{\beta}(x)\| < 1$ for all $x \in B(x^*,\alpha)$ and $\beta$ large enough, proving uniqueness of the equilibrium as well as assertion $(ii).$ The last assertion follows from Propositions \ref{prop:equilibria} and \ref{reversequil}.
\qed

A Nash equilibrium $x*$ is called {\em fully mixed } if $Supp(x^*) = \{1, \ldots, n \}$  and {\em partially mixed } if $ 1 < card(Supp(x^*)) < n.$

A fully mixed Nash equilibrium is called {\em non degenerate} if for all $u \in T \Delta$
$$ \left[ \forall w   \in T \Delta  \: \langle DU(x^*) u, w \rangle = 0 \right ] \Rightarrow  u = 0. $$
Let $$T \Delta (x^*) = \{u \in T \Delta: \: u_i = 0 \mbox{ for } i \not \in Supp(x^*)\}.$$
A partially mixed equilibrium $x^*$ is called non degenerate if
for all $u \in T \Delta(x^*)$
$$ \left[ \forall w   \in T \Delta(x^*)  \: \langle DU(x^*) u, w \rangle = 0 \right ] \Rightarrow  u = 0, $$
\blem
\label{lem:extrin} Let $x^* \in \Delta$ be a  mixed equilibria. Assume that $Supp(x^*) = \{1, \ldots, r\}$ for some $1 < r \leq n$ and set  $$x^* = (q_1,\ldots,q_{r-1},1-\sum_{i = 1}^{r-1} q_i, 0, \ldots,0).$$
Let, for $i = 1, \ldots, r-1,$
 \begin{eqnarray*}
 h^r_i(x_1, \ldots, x_{r-1},y_1, \ldots, y_{n-r}) &=& U_i(x_1, \ldots, x_{r-1}, 1-\sum_{i = 1}^{r-1} x_i - \sum_{i = 1}^{n-r} y_i, y_1, \ldots, y_{n-r}) \\
 & - &   U_r(x_1, \ldots, x_{r-1}, 1-\sum_{i = 1}^{r-1} x_i - \sum_{i = 1}^{n-r} y_i, y_1, \ldots, y_{n-r})
 \end{eqnarray*}
 Then $x^*$ is non degenerate if and only if the $(r-1) \times (r-1)$ matrix $$\left[\frac{\partial h^r_i}{\partial x_j}((q,0))\right]_{i,j = 1, \ldots r-1}$$ is invertible.
\elem
\prf One has $$\frac{\partial h^r_i}{\partial x_j}(q,0) = (\frac{\partial U_i}{\partial x_j}(x^*) - \frac{\partial U_i}{\partial x_r}(x^*))
- (\frac{\partial U_r}{\partial x_j}(x^*) -\frac{\partial U_r}{\partial x_r}(x^*)).$$
Let $$v = (v_1, \ldots, v_{r-1}, -\sum_{i = 1}^{r-1} v_i, 0, \ldots, 0) \in T\Delta(x^*)$$ and $$w = (w_1, \ldots, w_{r-1}, -\sum_{i = 1}^{r-1} w_i, 0, \ldots, 0) \in T\Delta(x^*).$$ Then it is easily seen that
$$\sum_{i = 1}^{r-1} \sum_{j=1}^{r-1} \frac{\partial h^r_i}{\partial x_j}(q) v_i w_j = \langle DU(x^*) v, w \rangle.$$ This proves that
$x^*$ is non degenerate if and only if $\left[\frac{\partial h^r_i}{\partial x_j}((q,0))\right]_{i,j = 1, \ldots r-1}$ is invertible. \qed
\bthm
\label{th:equimix}
Let $x^*$ be a non degenerate fully mixed Nash equilibrium for $U$ and  ${\cal N}$ a (sufficiently small)  neighborhood of $x^*.$  Then, there exists $\beta_0 > 0$ such that for all $\beta \geq \beta_0$ and $F \in \chi(\beta,U)$
$$\mathbf{Eq}(F) \cap {\cal N} = \{x^*_{\beta}\}.$$
Assume furthermore that the game is a potential game with potential $W.$ Then $x^*_{\beta}$ is hyperbolic for $F$ and its unstable  manifold (for $F$)  has dimension $\mathbf{Ind}(x^*,W|_{\Delta})$ under   one of the following conditions :
\bdes
\ita $L$  (hence $F$) is $C^2$ and $x^*_{\beta}$ is hyperbolic for $F,$  or
\itb  $F \in \chi_{rev}(\beta, U).$ \edes
\ethm
\prf
Set $T = 1/\beta.$
 Equilibria of $F_{\pi_{\beta}}$ are given by the set of equations
$$T (\log(x_i) - \log(x_n) ) = U_i(x) - U_n(x), i = 1, \ldots, n-1$$ or, with the notation of Lemma \ref{lem:extrin},
\beq
\label{eq:implic}
T(\log(x_i) - \log(1-\sum_{i = 1}^{n-1} x_i)) = h_i^n(x_1, \ldots, x_{n-1}), i = 1, \ldots, n-1.
\eeq
Write $x^* = (q_1, \ldots, q_{n-1}, 1-\sum_{i = 1}^{n-1} q_i).$
For $T = 0,$ $q = (q_1, \ldots, q_{n-1})$ is solution to (\ref{eq:implic}). Hence, by the implicit function theorem (which hypothesis is fulfilled  by the non degeneracy of $x^*$ and Lemma \ref{lem:extrin}) there exists $\alpha_0 > 0,$ a neighborhood $O$ of $q$ in $(\RR_+^*)^{n-1}$ and a $C^1$ map $T \in ]-\alpha_0,\alpha_0[ \mapsto q(T) \in O$ such that $(T, q(T))$ is the unique solution to (\ref{eq:implic}) in $]-\alpha_0,\alpha_0[ \times O.$
This proves the first assertion of the theorem with $\beta_0 > 1/\alpha_0$ and  $x^*_{\beta} = (q(1/\beta), 1-\sum_{i = 1}^{n-1} q_i(1/\beta))$

In case, the game is a potential game with potential $W,$ $F$ is gradient-like with Lyapounov function $V_{\beta}$ given by (\ref{eq:lyapgibbs}). Since $x^*$ is fully mixed, $\frac{1}{q_i} < \infty$ so that $\| \frac{1}{\beta} Hess V_{\beta}(x^*_{\beta}) - Hess W(x^*)\| \rightarrow 0$ as $\beta \rightarrow \infty.$ In particular, for $\beta$ large enough  $Hess V_{\beta}(x^*_{\beta})$ is non degenerate, because $x^*$ is non degenerate.
The last assertion then follows from  Propositions \ref{prop:equilibria} and \ref{reversequil}.
\qed
\bthm Let $x^*$ be a strict and non degenerate partially mixed Nash equilibrium for $U$ which support has cardinal $1 < r < n.$ Let ${\cal N}$ be  a (sufficiently small) neighborhood of $x^*.$ Then, there exists  $\beta_0 > 0$ such that for all $\beta \geq \beta_0$ and $F \in \chi(\beta,U)$
$$\mathbf{Eq}(F) \cap {\cal N} = \{x^*_{\beta}\}.$$
Assume furthermore that the game is a potential game with potential $W$ and that
 one of the following conditions hold :
\bdes
\ita $L$ (hence $F$)  is $C^2$ and $x^*_{\beta}$ is hyperbolic for $F,$  or
\itb  $F \in \chi_{rev}(\beta, U).$
\edes
Then $x^*_{\beta}$ is hyperbolic and $$k \leq \dim(E^u_{x^*_{\beta}}) \leq \min(n-r + k, r-1).$$
with $k = \mathbf{Ind}(x^*, W|_{\Delta(x^*)})$ and  $\dim(E^u_{x^*_{\beta}})$ stands for the dimension of the unstable manifold (for $F$).
\ethm
\prf
Assume without loss of generality that $Supp(x^*) = \{1,\ldots,r\}$ and set $x^* = (q_1,\ldots,q_{r-1},1-\sum_{i = 1}^{r-1} q_i, 0, \ldots, 0).$
 Write every element of $\Delta$ as $(x_1,\ldots, x_{r-1}, 1-\sum_{i = 1}^{r-1} x_i - \sum_{i = 1}^{n-r} y_i, y_1, \ldots y_{n-r})$ and  set $x = (x_1,\ldots,x_{r-1}), y = (y_1,\ldots,y_{n-r}).$
 Thus, with $\beta = 1/T,$ equilibria of $F_{\pi_{\beta}}$ are given by the following system of equations:

 \beq
  \label{eq:implimixt}
   T(\log(x_i) - \log(1-\sum_{i = 1}^{r-1} x_i - \sum_{i = 1}^{n-r} y_i)) = h^r_i(x,y), i = 1, \ldots r-1
 \eeq
 and
 \beq
 \label{eq:implicmixt2}
   T (\log(y_i) - \log(1-\sum_{i = 1}^{r-1} x_i - \sum_{i = 1}^{n-r} y_i)) = h^r_{i+r}(x,y), i = 1 \ldots n-r
 \eeq
 where $h^r_i$ is defined in Lemma \ref{lem:extrin}. The triplet $(T = 0,  x = q, y = 0)$ is solution to (\ref{eq:implimixt}).  Thus by the non degeneracy hypothesis and the implicit function theorem, there exists a smooth map $$\hat{x} : {\cal O} \mapsto {\cal V}, (T,y) \mapsto \hat{x}(T,y)$$ where ${\cal O}$ is a  neighborhood of $(0,0)$ in  $\RR \times \RR^{n-r}$ and ${\cal V}$ a neighborhood of $q$ in $\RR^{r-1}$ such that $(T, \hat{x}(T,y), y)$ is solution to (\ref{eq:implimixt}). Recall that $0 < \sum_{i = 1}^{r-1} q_i < 1$ and $h^r_{i+r}(q,0) < 0$ for all $i = 1, \ldots, n-r$ (because $x^*$ is strict). Thus,  by choosing ${\cal O}$ small enough we can furthermore ensure that
 \beq
 \label{hr1}
 0 < 1-\sum_{i = 1}^{r-1} \hat{x}_i(T,y) - \sum_{i = 1}^{n-r} y_i < 1
 \eeq
 and
 \beq
 \label{hr2}
 h^r_{i+r}(\hat{x}(T,y),y) \leq - \delta < 0, i = 1 \ldots  n-r
 \eeq
  for all $(T,y) \in {\cal O}.$

 Now replacing $x$ by $\hat{x}(T,y)$ in (\ref{eq:implicmixt2}) leads to
 $$y_i = G_i(T,y), i = 1 \ldots n-r$$
 where
 $$ G_i(T,y) =  (1-\sum_{i = 1}^{r-1} \hat{x}_i(T,y) - \sum_{i = 1}^{n-r} y_i) \exp{(\frac{1}{T} h^r_{i+r} (\hat{x}(T,y),y))}.$$
Using (\ref{hr1}) and (\ref{hr2}) we see that $\alpha$ small enough and $T \leq \frac{\log(1/\alpha)}{\delta}$   $G(T, \cdot)$ maps $\{y  \in \RR^{n-r}: \: 0 \leq y_i \leq \alpha\}$ into itself. By Brouwer's fixed point theorem, $G(T,\cdot)$ admits a fixed point $\hat{y}(T).$ Furthermore, $\|D_y G(T,y)\| \leq \frac{C}{T} e^{-\delta/T}$ for some constant $C $ making $G(T,\cdot)$ a contraction. This implies that $\hat{y}(T)$ is unique. Finally define $x^*_{\beta}$ by  $x^*_{\beta,i} = \hat{x}_i(T,\hat{y}(T))$ for $1 \leq i < r$ and $x^*_{\beta,i+r} =  \hat{y}_i(T)$ for $1 \leq 1 \leq n-r.$ \\

 We now prove the last assertions.  By assumption, $T \Delta(x^*)$ admits a decomposition
  $T \Delta(x^*) = E_{+} \oplus E_{-}$ with $\langle Hess(W)(x^*) u, u \rangle > 0$  (respectively $ < 0 $) for all $u \in E_+$ (respectively $E_-$) and $u \neq 0.$

Set $T \Delta_s(x^*) = \{u  \in T\Delta : \: u_1 = \ldots = u_r = 0\}.$
Then $$T\Delta = E_{+} \oplus E_{-} \oplus T \Delta_s(x^*).$$
 Let now $V_{\beta}$ be the Lyapounov function given (\ref{eq:lyapgibbs}). Then  for all $u \in T\Delta$
 $$Q_{\beta}(u) := \langle \frac{1}{\beta} Hess(V_{\beta})(x^*_{\beta}) u, u \rangle = \langle Hess W(x^*_{\beta}) u, u \rangle  + \frac{1}{\beta} \sum_i \frac{1}{x^*_{\beta,i}} u_i^2.$$
The construction of $x^*_{\beta}$ shows that $\frac{1}{\beta} \frac{1}{ x^*_{\beta,i}} \rightarrow 0$ for  $i \leq r$ and  $\frac{1}{\beta}
\frac{1}{ x^*_{\beta,i}}  \rightarrow \infty$ for $i > r$ when $\beta \rightarrow \infty.$
 Thus, for $\beta$ large enough,
$Q_{\beta}$ is non degenerate,   definite positive on $E_+$ and $T\Delta_s(x^*),$ and  definite negative on $E_-.$

This implies that its index is bounded below by $k = \dim(E_{-})$ and above by $\min{(r-1, n-r - k)}.$  This index equals the dimension of the stable manifold by Proposition \ref{prop:equilibria}. Under the reversibility assumption hyperbolicity follows from  Proposition \ref{reversequil}.
\qed
\section{Reversibility and Gradient Structure}
\label{sec:reversibility}
Recall that an irreducible rate matrix $L$ is called {\em reversible} with respect to $\pi \in \dot{\Delta}$ is $\pi_i L_{ij} = \pi_j L_{ji}.$
In this case $\pi$ is the (unique) invariant probability of $L.$ Here we will consider gradient properties of (\ref{ode}) under the assumption that $L(x)$ is reversible.

A $C^k, k \geq 0$ (Riemannian) metric on $\dot{\Delta}$ (or $\Delta$) is a $C^k$ map $g$ such that for each $x \in \Delta$
$g(x) : T\Delta \times T \Delta \mapsto \RR$ is a definite positive bilinear form.
Given a $C^1$ map $V : \dot{\Delta} \mapsto \RR$ we let $grad_g V$ denote the gradient vector field of $V$ with respect to $g.$
That is
$$g(x)(grad_g V(x),u) = \langle \nabla V(x) , u \rangle$$ for all $u \in T\Delta.$

\bprop
Assume that for all $x \in \dot{\Delta}$ $L(x)$ is reversible with respect to $\pi(x)$ and assume that the map $h : \dot{\Delta} \mapsto \RR^n,$ defined by $$h(x) = \frac{x}{\pi(x)}$$ is a $\alpha$-quasigradient. Then there exists a  metric  $g$ on $\dot{\Delta}$ such that  for all $x \in \dot{\Delta}$ $F(x) = - grad_g V(x).$ If $L$ and $\alpha$  are $C^k$ then $g$ is $C^k.$
\eprop
\prf
The proof is similar to the proof of Proposition \ref{reversequil}.
 Let $A(x) : T \Delta \mapsto T \Delta$ be defined by $A(x) h = - h L(x).$ Then $A(x)$ and $L(x)$ are conjugate by the relation
$\pi(x) L(x) h = A(x) \pi(x) h$ and $A(x)$ is a definite positive operator for  the scalar product on $T \Delta$ defined by
$\langle u, v \rangle_{1/\pi(x)} = \sum_i u_i v_i \frac{1}{\pi_i(x)}.$ Define now a Riemannian metric on $T\Delta$ by
\beq
\label{eq:metric} g_0(x)(u,v) = \langle A(x)^{-1} u, v \rangle_{\frac{1}{\pi(x)}}.
\eeq
Since
$F(x) = x L(x) = (x -\pi(x)) L(x) = A(x)(-x + \pi(x)),$ we get
$$g_0(x)(F(x),u) = - \langle \frac{x}{\pi(x)} -1, u \rangle = - \langle \frac{x}{\pi(x)} , u \rangle.$$
If $x \mapsto \frac{x}{\pi(x)}$ is a quasi gradient, this makes $F$ a gradient for the metric $g(x) = \alpha(x) g_0(x).$
\qed
\bex {\rm Suppose that  $L(x)$ is reversible with respect to $\pi,$ independent on $x.$ Then $x \mapsto \frac{x}{\pi}$ is the gradient of the $\chi^2$ function $V(x) = \sum_i (\frac{x_i}{\pi_i} - 1)^2 \pi_i.$ Hence $F(x) = - grad_g V(x)$  for some metric $g$.}
\eex

Under the weaker assumption that $x \mapsto s(\frac{x}{\pi(x)})$ is a quasi-gradient for some strictly increasing function $s$ (see Theorem \ref{main})) it is no longer true that $F$ is a gradient, but it can be approximated by a gradient.  The next Lemma is the key tool. Its proof is identical to the proof of Proposition \ref{reversequil}.
\blem
\label{reversequil2}
Assume that assumptions of Theorem \ref{main} hold and that
 for all $x \in \dot{\Delta}, L(x)$ is reversible with respect to $\pi(x).$ Then there exists a  metric $g_0$ on $\dot{\Delta}$ such that for $p \in \mathbf{Eq}(F) \cap \dot{\Delta}$ and  $u,v \in T\Delta$
$$g_0(p)(DF(p) u, v) = - \langle Hess(V)(p) u, v \rangle.$$ If, furthermore, $L$ and $\alpha$ (in equation \ref{defV}) are $C^k$ then $g_0$ is $C^k.$
\elem

\bthm
\label{th:closetograd}
Assume that
\bdes
\ita  Assumptions of Theorem \ref{main} hold with $s, \alpha$ and $L$  $C^k, k \geq 2,$
\itb  For all $x \in \dot{\Delta}$
$L(x)$ is reversible with respect to $\pi(x),$
\itc $\mathbf{Eq}(F) \cap \dot{\Delta}$ is finite.
\edes
Then for every neighborhood ${\cal U}$ of $\mathbf{Eq}(F) \cap \dot{\Delta}$ and every $\eps > 0$ there exists a $C^k$ metric $g$ on $\dot{\Delta}$
such that
\bdes
\iti $- grad_g V = F$ on $\dot{\Delta} \setminus {\cal U}.$
\itii $\|- grad_g V -F \|_{C^1 , {\cal U}} \leq \eps$
where $$\| G \|_{C^1 , {\cal U}} = \sup_{x \in  {\cal U} } \|G(x)\| + \|DG(x)\|.$$
\edes
\ethm
\prf
Let ${\cal E} = \dot{\Delta} \cap \mathbf{Eq}(F),$
$v(x) = d(x,{\cal E} ) = \min_{ p \in {\cal E}} \|x-p\|$
 and let $\psi : \RR^+ \mapsto [0,1]$ be a $C^{\infty}$ function which is $0$ on $[0,1],$  $1$ on $[3, \infty[$ and such that $0 \leq \psi' \leq 1.$ Fix $\eps > 0$ and let $\lambda(x) = \psi(\frac{v(x)}{\eps}),$
  $G_0 = - grad_{g_0} V$ where $g_0$ is given by Lemma \ref{reversequil2} and   $$G(x) = (1-\lambda(x)) G_0(x) + \lambda(x) F(x).$$ Since for all $p \in {\cal E}, F(p) - G_0(p) = DF(p) -DG_0(p) = 0$ there exists a constant $C > 0$ such that
  $$\|G_0(x) - F(x)\| \leq C v(x)^2, \; \|DG_0(x) - DF(x)\| \leq C v(x).$$ Thus
 $$\|G(x) - F(x)\| = (1-\lambda(x))\|G_0(x)-F(x)\| \leq C (1-\lambda(x)) v(x)^2 \leq C \eps^2$$ and
 $$\|DG(x) - DF(x)\| = \|(1-\lambda(x)) (DG_0(x) -DF(x)) + \langle \nabla \lambda(x), G_0(x) - F(x) \rangle \| $$
 $$\leq C ((1-\lambda(x)) v(x) + \frac{1}{\eps} v(x)^2) \leq C \eps.$$
This shows that $G$ is a $C^1$ approximation of $F$ which coincides with $F$ on $\{ v(x) \geq 3 \eps\} $ and with $G_0$ on $\{v(x) \leq \eps\}.$
Furthermore, $$\langle \nabla V(x), G(x) \rangle =  - (1-\lambda(x)) g_0(x)(G_0(x),G_0(x)) + \lambda(x) \langle \nabla V (x), F(x) \rangle \leq 0$$ with equality if and only if $x \in {\cal E}.$

Now, for all $x \in \dot{\Delta} \setminus {\cal E}$ $$T\Delta = \nabla V(x)^{\bot} \oplus \RR G(x)$$ and the splitting is smooth in $x.$
Hence $u \in T\Delta$ can be uniquely written as $u =  P_x(u) + t_x(u) G(x)$ with $t_x(u) \in \RR$ and $P_x(u) \in \nabla V(x)^{\bot}.$
Let  $g$ be the metric on $\dot{\Delta} \setminus {\cal E}$ defined by
$$g(x)(u,v) = g_0(P_x(u),P_x(v)) + t_x(u)t_x(v) g_0(x)(G_0(x),G(x)).$$ Then $g$ coincides with $g_0$ on $\{0 < x < v(x) < \eps\}$ so that
 $g$ can be extended to a $C^2$ metric on $\dot{\Delta}.$ By construction of $G$ and $g,$  $G = -grad_g V.$
 \qed

\section{Questions of Structural Stability}
\label{sec:struct}
Let $C^k_{pos}(\Delta,T\Delta)$ denote the set of $C^k$ vector fields
  $F : \Delta \mapsto T\Delta$ leaving $\Delta$ positively invariant.

Two elements $F,G \in C^k_{pos}(\Delta,T\Delta)$  are said
{\em topologically equivalent} if there exists a homeomorphism $h: \Delta \mapsto \Delta$ which takes orbits of
$F$ to orbits of $G$ preserving their orientation. A set $\chi \subset C^k_{pos}(\Delta,T\Delta)$ is said {\em structurally stable} if all its elements are topologically equivalents.

Let $\pi : \Delta \mapsto \dot{\Delta}$ be a smooth function. Assume that $\pi$  verifies the assumption of Theorem  \ref{main} and  that $F_{\pi}$ has non degenerate equilibria. Let  $\chi_{\pi, \, rev}$ denote the convex set of vector fields having the form given by (the right hand side of) (\ref{ode}), where for each $x \in \Delta,$ $L(x)$ is irreducible and reversible with respect to $\pi(x).$
 By Theorem \ref{main}, Proposition \ref{reversequil} and Theorem \ref{th:closetograd} all the elements of $\chi_{\pi, \, rev}$  have the same strict Lyapounov function $V,$  hyperbolic equilibria (given by the critical points of $V$) and are $C^1$ close to $-grad_g V$ for some metric $g.$ We may then wonder wether  $\chi_{\pi, \, rev}$  is structurally stable. The following construction shows that this is not the case.

 \subsection{Potential Games are not structurally stable}
 Here $\Delta$ stands for the $2$-dimensional simplex in $\RR^3.$
Let $$\tilde{\Delta} = \{(y_1,y_2) \in \RR^2 \: : y_1, y_2 \geq 0, \: y_1 + y_2 \leq 1\}$$
and
 $\jmath :  \RR^3  \mapsto  \RR^2$ be the projection defined by $\jmath(x_1,x_2,x_3) = (x_1,x_2).$ Note that $\jmath$ maps $\Delta$ homeomorphically onto $\tilde{\Delta}.$

Let $\tilde{W} : \RR^2 \mapsto \RR$ be a smooth function.
 Assume that
 \bdes
 \ita  $- \nabla \tilde{W}$ points inward $\tilde{\Delta}$ on $\partial \tilde{\Delta};$
 \itb The critical set  $crit(\tilde{W}) = \{y  \in \tilde{\Delta} :\: \nabla \tilde{W}(y) = 0 \}$ consist of (finitely many)  non degenerate points,
 \itc For all $u \in \RR$ $$\frac{\partial \tilde{W}}{\partial y_1}(u,u) = \frac{\partial \tilde{W}}{\partial y_2}(u,u).$$ In particular, the diagonal $D(\tilde{\Delta}) = \{ (y_1, y_2) \in \tilde{\Delta}: \: y_1 = y_2\}$ is positively invariant under the dynamics \beq \label{eq:odetW}
 \dot{y} = - \nabla \tilde{W}(y) \eeq
 \itc There is a saddle connection contained in $D(\tilde{\Delta}),$ meaning that there are two saddle points of $\tilde{W}$ $s^1, s^2 \in D(\tilde{\Delta})$ and some (hence every) point $y \in ]s^1,s^2[$ which $\alpha$ limit set under (\ref{eq:odetW}) is $s^1$ and omega limit set is $s^2.$
 \edes
 It is not hard to construct such a map.\\

 Let $W : \RR^3 \mapsto \RR$ be defined by $W = \tilde{W} \circ \jmath.$

 Consider now the $3-$ strategies potential game associated to $W.$
 Payoffs are then defined by $$U_i(x) = - \frac{\partial \tilde{W}}{\partial x_i}(x_1,x_2), i = 1, 2 \mbox{ and }
 U_3(x) = 0.$$
Using  the notation of section \ref{sec:equilibriapot}, record that  $F_{\pi_{\beta}} = -Id + \pi_{\beta}$ where $\pi_{\beta}$ is defined by (\ref{defpibeta}), and $\chi_{rev}(\beta,U)$ is the set of vector fields given by (\ref{ode}) with $L(x)$ irreducible and reversible with respect to $\pi_{\beta}(x).$
\bprop
\label{prop:counter}
For all $\beta > 0$ sufficiently large, there exists $F \in \chi_{rev}(\beta,U)$ (which can be chosen $C^1$ close to $F_{\pi_{\beta}}$ ) which is not equivalent to $F_{\pi_{\beta}}.$
\eprop
\subsubsection*{Proof of Proposition \ref{prop:counter}}
By definition of Nash equilibria (see section \ref{sec:equilibriapot}) and condition $(a)$ above, Nash equilibria of $U$ are fully mixed and coincide with critical points of $\tilde{W}:$
 $$crit(\tilde{W}) =  \jmath(\mathbf{NE}(U)).$$
  \blem
  \label{lem:contrex}
  For all $\eps > 0$ there exists $\beta_0 > 0$ such that for all  $\beta \geq \beta_0$ and  $F \in \chi_{rev}(\beta,U)$ there is a one to one map
  $$p \in crit (\tilde{W})  \mapsto p_{\beta} \in \mathbf{Eq}(F),$$
  such that
  \bdes
   \iti $\|p-\jmath(p_{\beta})\| \leq \eps,$
   \itii The unstable (respectively stable) manifold of $p_{\beta}$ has dimension $\mathbf{Ind}(p,\tilde{W},)$ (resp. $2-\mathbf{Ind}(p,\tilde{W})$). In particular, $s^1_{\beta}$ and $s^2_{\beta}$ are saddle points.
   \itiii   $p \in  D(\tilde{\Delta})  \Leftrightarrow \jmath(p_{\beta}) \in D(\tilde{\Delta})$
   \itiv Under the dynamics induced by $F_{\pi_{\beta}},$ the interval $[s^1_{\beta}, s^2_{\beta}]$ is invariant and for some (hence all) $q \in ]s^1_{\beta},  s^2_{\beta}[$ the alpha limit (respectively omega limit) set of $q$ equals $s^1_{\beta}$ (respectively $s^2_{\beta}$).
  \edes
  \elem
  \prf Assertions $(i)$ and $(ii)$ this follows from Propositions  \ref{prop:equiNash} and \ref{th:equimix}.

  On $\jmath^{-1}(D(\tilde{\Delta})) = \{(x_1,x_1, 1-2x_1)\}$ equilibria of $F_{\pi_{\beta}}$ are given by the implicit equation
  $T (\log(x_1) - \log(1-2 x_1)) = U_1(x_1,x_1)$ where $T = 1/\beta.$ Solutions for $T = 0$ coincide with $\jmath^{-1} (D(\tilde{\Delta}) \cap crit(\tilde{W})).$ For $T > 0$ and small enough, assertion $(iii)$ then follows from the  implicit function theorem.

  By condition $(c),$ $\frac{\partial  \tilde{W}}{\partial x_1} = \frac{\partial  \tilde{W}}{\partial x_2}$ on $D(\tilde{\Delta}).$ Thus $U_1(x) = U_2(x)$ (hence $F_{\pi_{\beta,1}}(x) = F_{\pi_{\beta,2}}(x))$ on $\jmath^{-1}(D(\tilde{\Delta}))$ proving invariance of $[s^1_{\beta} , s^2_{\beta} ] \subset \jmath^{-1}(D(\tilde{\Delta})).$ Assertion $(iv)$ follows since, by $(iii),$ there are no equilibria in $]s^1_{\beta}, s^2_{\beta}[.$
  \qed

We now construct $F \in \chi_{rev}(\beta,U).$
Let $L(x)$ be the rate matrix defined for $i \neq j$ by $$L_{ij}(x) = \pi_{\beta,j}(x) \mbox{ if } i,j \not \in \{1,3\}$$
$$L_{13}(x) = (1 + a(x)) \pi_{\beta,3}(x) \mbox{ and } L_{31}(x) = (1+a(x)) \pi_{\beta,1}(x)$$
where $a : \Delta \mapsto \RR^+$ is a smooth function to be defined below. Then
equation (\ref{ode}) reads
\begin{eqnarray*} \large
  \dot{x_1} &=& (x_2 \pi_{\beta,1}(x) - x_1  \pi_{\beta,2}(x))  + (x_3 \pi_{\beta,1}(x) - x_1  \pi_{\beta,3}(x))(1+ a(x)), \\
  \dot{x_2} &=& (x_1 \pi_{\beta,2}(x) - x_2  \pi_{\beta,1}(x))  + (x_3 \pi_{\beta,2}(x) - x_2  \pi_{\beta,3}(x)),\\
  \dot{x_3} &=& -\dot{x_1} - \dot{x_2}.
\end{eqnarray*}
Thus, on $x_1 = x_2,$
$$\dot{x_1} - \dot{x_2} =  [x_3 \pi_{\beta,1}(x) - x_1 \pi_{\beta,3}(x)] a(x)$$
  $$= \frac{a(x)}{Z(x)} (x_3 e^{\beta U_1(x)} - x_1).$$ The map $x \mapsto x_3 e^{\beta U_1(x)} - x_1$ vanishes at points $s_{\beta}^1, s_{\beta}^2$ and has a constant sign over $[s^1_{\beta}, s^2_{\beta}]$ (for otherwise there would exists an equilibrium for $F$ in $]s^1_{\beta}, s^2_{\beta}[$ contradicting Lemma \ref{lem:contrex}).
Let $p = (s^1_{\beta}+ s^2_{\beta})/2$ and $B_{\eta}$ be the Euclidean open ball with center $p$ and radius $\eta.$ Choose $\eta$ small enough so that
\bdes
\iti $B_{\eta} \cap [s^1_{\beta}, s^2_{\beta}] = ]q^1,q^2[$ with $s^1_{\beta} < q^1 < q^2 < s^2_{\beta}$ where $<$ stands for the natural ordering on $[s^1_{\beta}, s^2_{\beta}].$
\itii $x \mapsto x_3 \pi_{\beta,1}(x) - x_1 \pi_{\beta,3}(x)$ has constant sign on $B_{\eta}.$
\edes
Let $x \mapsto a(x)$ be such $a = 0$ on $\Delta \setminus B_{\eta},$  $a > 0$ on $B_{\eta}$ and $0 \leq a \leq \eta$ on $\Delta.$ Then,
 the alpha limit set of $q^1$ equals $s^1_{\beta},$ for both $F$ and $F_{\pi_{\beta}}$ but since $\dot{x_1} - \dot{x_2}$ doesn't vanish on $B_{\eta}$ the trajectory through $q^1$ exits $B_{\eta}$ at a point $\neq q^2$ and, consequently,  the omega limit set of $q^1$ for $F$ is distinct from $s^2_{\beta}.$ This proves that $F$ and $F_{\pi_{\beta}}$ are not equivalent.
 \subsection{Open question}
 The preceding construction shows that $\chi_{rev}( \beta,U)$ is not structurally stable for an arbitrary potential game but this might be the case for particular examples. Consider for example the  Gibbs model described in Remark \ref{rem:DF}.
 For $U^0 \in \RR^n$ and $U = (U_{ij})$ symmetric, let $\chi_{rev}(\beta, U^0, U)$ be the set of $C^1$ vector field given by (\ref{ode}) with $L(x)$ irreducible and reversible with respect to the Gibbs measure (\ref{Gibbs}).
 \paragraph{Question} For generic $(U^0,U)$ and $\beta$ large enough, is $\chi_{rev}(\beta, U^0, U)$ structurally stable ?

\section{Appendix}
Let $L$ be an irreducible rate matrix and $\pi \in \dot{\Delta}$ denote the invariant probability of $L$. That is the unique solution (in $\Delta$)
of $\pi L = 0.$
For all $f,g \in    \mathbb{R}^n$  we let $$\langle f, g \rangle = \sum_i f_i g_i, \, \langle f, g \rangle_{\pi} = \sum_i f_i g_i \pi_i \mbox{ and } \langle f, g \rangle_{1/\pi} = \sum_i f_i g_i {\frac{1}{\pi_i}}.$$ The {\em Dirichlet form} of $L$ is the map ${\cal E} :  \mathbb{R}^n \mapsto  \mathbb{R}_{+}$ defined as
   $${\cal E}(f) = - \langle f, L f\rangle_{\pi} =  \frac{1}{2} \sum_{i,j} (f_i-f_j)^2 L_{ij} \pi_i.$$
   By irreducibility, ${\cal E}(f) > 0$ unless $f$ is constant, and the {\em spectral gap}
   $$\lambda = \sup \{ {\cal E}(f)  : \: \langle f, 1 \rangle_{\pi} = 0, \langle f, f \rangle_{\pi} = 1 \}$$ is positive.
   We let $L^*$ be the irreducible rate matrix defined by
   $$L^*_{ij} = \frac{\pi_j L_{ji}}{\pi_i}.$$ Note that $L^*$ admits $\pi$ as invariant probability and that $L^*$ is the adjoint of $L$ for $\langle, \rangle_{\pi}.$

   We let $L^T : T\Delta \mapsto T \Delta$ be defined by $$L^T h = h L.$$
   Finally recall that for all $f \in \RR^n$ $\frac{f}{\pi}$ stands for the vector defined by $(\frac{f}{\pi})_i = \frac{f_i}{\pi_i}, i = 1 \ldots n.$
   \blem
   \label{lemonLT}
   For all $u,v \in T \Delta$
   $$\langle  L^T u, v \rangle_{1/{\pi}}  = \langle L^* (\frac{u}{\pi}), \frac{v}{\pi}) \rangle_{\pi}$$
   In particular $L^T$ is invertible and $L^T$ is a definite negative operator for $\langle , \rangle_{\frac{1}{\pi}}$ whenever $L$ is reversible with respect to $\pi.$
   \elem
   \prf The first assertion follows from elementary algebra. For the second, note that $\langle  L^T u, u \rangle_{1/{\pi}} = - {\cal E}(\frac{u}{\pi}).$ Thus, by irreducibility, $$\langle  L^T u, u \rangle_{1/{\pi}} < 0$$ unless $u = 0.$

   \qed

\subsection{Proof of Lemma \ref{poincare1}}
Given $f \in  \mathbb{R}^n$    we write $f \geq 0$ if $f_i \geq 0$ for all $i.$    We let $1 \in     \mathbb{R}^n$ denote the vector
which components are all equal to $1.$
 For all $t \geq 0$ we let $P_t = e^{t L}.$ Since $L$ is a rate matrix, $(P_t)$ is a Markov semigroup meaning that
$P_t f \geq 0$   for all $f \in \mathbb{R}^n$ with $f \geq 0$ and $P_t 1 = 1.$

 \blem
   \label{lem1}
   Let $I \subset \mathbb{R}$ be an open interval and $S : I \mapsto \mathbb{R}$ a $C^2$ function such that $S''(t) \geq \alpha > 0.$
   Let $f \in \mathbb{R}^n$ be such that $f_i \in I$ for all $i.$ Then
    $$\frac{d}{dt} \langle S(P_tf), 1\rangle_{\pi}|_{t = 0} \leq -\alpha {\cal E} (f).$$
   \elem
   \prf     For all $u,v \in I$ $S(v) - S(u) - S'(u)(v-u) \geq \alpha/2 (v-u)^2.$
    Hence for all $i,j$
    $$S(f_j) - S((P_tf)_i)   - S'((P_t f)_i) (f_j - (P_t f)_i) \geq \alpha/2 (f_j - (P_t f)_i)^2.$$
    Applying $P_t$ to this inequality gives
    $$P_t(Sf)_i - S((P_t f)_i) \geq \alpha/2 P_t (f_i-(P_tf)_i)^2 ) = \alpha/2 (P_tf^2_i- (P_tf)_i^2)$$
    Hence
        $$P_t(Sf) - S((P_t f)) \geq \alpha/2 P_t (f-(P_tf))^2 ) = \alpha/2 (P_tf^2- (P_tf)^2).$$
        Therefore, using the fact that $\langle P_t g, 1 \rangle_{\pi} = \langle g,1\rangle_{\pi}$ leads to
        $$ \langle Sf - S(P_tf), 1\rangle_{\pi} \geq \alpha \langle f^2 - (P_t f)^2, 1 \rangle_{\pi}.$$     Dividing by $t$ and letting $t \rightarrow 0$ leads     to the desired inequality. \qed

        Let $S : ]0,\infty[ \mapsto \RR$ be a $C^2$ function with positive second derivative.
 Let $H^S_\pi : \Delta \mapsto \mathbb{R}$ be the map defined by
$$H^S_{\pi}(x) = \sum_i \pi_i S(\frac{x_i}{\pi_i}).$$
\bcor
\label{cor1}
For all $x \in  \Delta$ $$\langle \nabla H_{\pi}^S(x),xL \rangle \leq -\alpha \lambda Var_{\pi}(f) $$
 where $f_i = \frac{x_i}{\pi_i}$
\ecor
\prf For $x \in \Delta$ let $x(t) = x e^{tL},$ $f_i = \frac{x_i}{\pi_i}, \, f_i(t) = \frac{x_i(t)}{\pi_i}$ and
$P^*_t g = e^{tL^*} g.$
Note that  $P^*_t$) is the adjoint of  $P_t$ with respect to $\langle, \rangle_{\pi}.$

For all $g \in \mathbb{R}^n,  \langle x(t), g \rangle = \langle x, P_t g \rangle = \langle f, P_t g \rangle_{\pi} = \langle P_t^* f , g \rangle_{\pi}$
so that $f(t) = P^*_t f.$ Hence by the preceding lemma applied to $L^*$ it follows that
     $$\langle \nabla H_{\pi}^S(x),xL \rangle =      \frac{d}{dt} \langle S(P^*_tf), 1\rangle_{\pi}|_{t = 0} \leq -\alpha {\cal E} (f) \leq -\alpha \lambda Var_{\pi}(f)$$  where $\alpha = \min_{i} S''(\frac{x_i}{\pi_i}) > 0.$
     \qed
We now prove the Lemma. Set $S(t) = \int_1^t s(u) du.$ Then for all $u \in T \Delta$
$$\langle \nabla H^S_{\pi}(x),u \rangle = \sum_i u_i s(\frac{x_i}{\pi_i})$$ and the results  follows from Corollary \ref{cor1}.

\bibliographystyle{amsplain}
\bibliographystyle{imsart-nameyear}
\bibliographystyle{nonumber}
\bibliography{lyapounov}

\section*{Acknowledgments}
This work was supported by the SNF grant $2000020_149871/1$
I would like to thank J. B Bardet, F. Malrieu, M. W Hirsch, J. Hofbauer, J. Robbin, B. Sandholm, S. Sorin, P. A Zitt for numerous discussions on topics
 related to this paper.
\end{document}